\documentclass[11pt,a4paper,twoside]{amsart}

\usepackage[utf8]{inputenc}
\usepackage[english]{babel}
\usepackage[T1]{fontenc}
\usepackage{microtype, fancyhdr, lmodern, xcolor}
\usepackage[a4paper, hmarginratio=1:1]{geometry}

\usepackage{amsfonts, amsmath, amsthm, amssymb, mathrsfs, amscd}		
\numberwithin{equation}{section}	
\usepackage{faktor}
\allowdisplaybreaks

\usepackage{url, enumitem}
\usepackage[noadjust]{cite}

\usepackage{multirow,bigdelim, caption}
\usepackage{booktabs}

\usepackage{graphicx}

\usepackage{hyperref}

\theoremstyle{plain}
\newtheorem{thm}{Theorem}[section]
\newtheorem{prop}[thm]{Proposition}
\newtheorem{cor}[thm]{Corollary}
\newtheorem{lem}[thm]{Lemma}

\theoremstyle{remark}
\newtheorem{rmk}[thm]{Remark}

\theoremstyle{definition}
\newtheorem{defn}[thm]{Definition}
\newtheorem{exam}[thm]{Example}

\newcommand\Zz{\mathbb{Z}}

\newcommand\Rr{\mathbb{R}}
\newcommand\Cc{\mathbb{C}}

\DeclareMathOperator{\Aut}{Aut}	

\newcommand\BB[1]{B_{#1}}	
\newcommand\PB[1]{P_{#1}}	
\newcommand\WB[1]{WB_{#1}}	
\newcommand\WBE[1]{WB_{#1}^{ext}}	
\newcommand\BP[1]{BP_{#1}}	
\newcommand\PC[1]{PC_{#1}}	
\newcommand\VB[1]{VB_{#1}}	
\newcommand\PLBE[1]{PLB_{#1}^{ext}}	
\newcommand\LBE[1]{LB_{#1}^{ext}}	
\newcommand\PLB[1]{PLB_{#1}}	
\newcommand\LB[1]{LB_{#1}}	
\newcommand\rB[1]{rB_{#1}}	
\newcommand\rBE[1]{rB_{#1}^{ext}}	

\newcommand \ep{\varepsilon}

\newcommand \F[1]{F_{#1}}

\newcommand\Gr[1]{\langle#1\rangle}		

\DeclareMathOperator{\Homeo}{Homeo}	
\DeclareMathOperator{\Diffeo}{Diffeo}	
\DeclareMathOperator{\PHomeo}{PHomeo}	
\DeclareMathOperator{\PDiffeo}{PDiffeo}	

\newcommand\id{\mathrm{id}}
\newcommand\ie{{\textit{i.e.}}}
\newcommand\ii{i}
\newcommand\inv{^{-1}}

\newcommand\jj{j}

\newcommand\MCG[2]{\mathrm{MCG}({#1}, {#2})} 	
\newcommand\PMCG[2]{\mathrm{PMCG}({#1}, {#2})} 	

\newcommand\nn{n}
\newcommand\nno{{n-1}}

\makeatletter
\newcommand{\oset}[2]{%
  {\mathop{#2}\limits^{\vbox to -.5\ex@{\kern-\tw@\ex@
   \hbox{\scriptsize #1}\vss}}}}
\makeatother

\newcommand\PR[1]{PR_{#1}}	
\newcommand\PUR[1]{PUR_{#1}}	
\newcommand\PRi[1]{\mathcal{PR}_{#1}}	
\newcommand\PUnRi[1]{\mathcal{PUR}_{#1}}	

\newcommand\R[1]{R_{#1}}	
\newcommand\Ri[1]{\mathcal{R}_{#1}}	

\newcommand\rr[1]{\rho_{\hspace{-0.3ex}#1}^{\null}}	
\newcommand\rrq[2]{\rho_{\hspace{-0.3ex}#1}^{#2}}	

\newcommand\sig[1]{\sigma_{\hspace{-0.3ex}#1}^{\null}} 
\newcommand\sigg[2]{\sigma_{\hspace{-0.3ex}#1}^{#2}}	
\newcommand\siginv[1]{\sigg{#1}{-1}}	

\newcommand\UR[1]{UR_{#1}}	
\newcommand\UnRi[1]{\mathcal{UR}_{#1}}	

\begin{document}

\title[A journey through loop braid groups]{A journey through loop braid groups}

\author[Damiani]{Celeste Damiani}
\address{Laboratoire de Math\'ematiques Nicolas Oresme, CNRS UMR 6139, Universit\'e de Caen BP 5186,  F-14032 Caen, France.
}
\email{celeste.damiani@unicaen.fr}



\date{\today}

\begin{abstract}
In this paper we introduce distinct approaches to \emph{loop braid groups}, a generalization of braid groups, and unify all the definitions that have appeared so far in the literature, with a complete proof of the equivalence of these definitions. These groups have in fact been an object of interest in different domains of mathematics and mathematical physics, and have been called, in addition to loop braid groups, with several names such as of \emph{motion groups}, groups of \emph{permutation-conjugacy automorphisms},  \emph{braid-permutation groups}, \emph{welded braid groups} and \emph{untwisted ring groups}. 
In parallel to this, we introduce an extension of these groups that appears to be a more natural generalization of braid groups from the topological point of view. 
Throughout the text we motivate the interest in studying loop braid groups and give references to some of their applications. 
\end{abstract}

\maketitle
\section{Introduction}
Braid groups $\BB\nn$ were introduced by Hurwitz \cite{Hur} in 1891 as fundamental groups of configuration spaces of $\nn$ points in the complex plane. However, they owe their name to Artin~\cite{Art}: he considered them in terms of braid automorphisms of~$\F\nn$, the free group of rank~$\nn$, but also in geometric terms. The geometric interpretation certainly is the most intuitive and best known, in particular because of its use in knot theory. Then, Magnus~\cite{Mag} considered braid groups from the point of view of mapping classes, while Markov~\cite{Mar} introduced these groups from a purely group-theoretic point of view. All these points of view have long been known to be equivalent~\cite{Zar}. We can then say that braid groups are ubiquitous objects. Any different definition carries a possible generalization; for instance, we can see braid groups as particular case of Artin-Tits groups, Garside groups, mapping class groups and surface braid groups. Few of these generalizations share with braid groups their principal property: a large family of different equivalent definitions.


\emph{Loop braid groups} $\LB\nn$ are a remarkable exception to this fact. Their study has been widely developed during the last twenty years.
 The first curious fact about these groups is that they appear in the literature with a large number of different names. The name \emph{loop braid groups} was proposed by  Baez-Crans-Wise~\cite{BWC}, while McCool and Savushkina~\cite{McC,SAV} considered them as \emph{basis-conjugating automorphisms} and \emph{permutation conjugacy automorphisms} of~$\F\nn$,  Fenn-Rim\'{a}nyi-Rourke~\cite{FRR} introduced the name of \emph{welded braids} (probably the most common notation nowadays), and Brendle-Hatcher~\cite{BH} called them \emph{groups of untwisted rings}. We will spend some words on the different appearances in Section~\ref{S:other}.

The first and principal aim of this paper is to introduce several interpretations of \emph{loop braid groups}, and provide a complete proof of the equivalence of these formulations. This provides a clear and well-established theoretical setting, as the one we have in the case of braid groups~$\BB\nn$. Although the equivalence of some of the formulations has already been proved or at least stated in the literature, for example in~\cite[Section 2, Theorem 2.6]{Bn}, we did not find the explicit isomorphisms between all of them.

Our second purpose is to introduce \emph{extended loop braid groups}~$\LBE\nn$, a generalization of loop braid groups. In the literature these groups have been less studied, but, as it will appear  throughout this work, they are even more interesting and rising in a more natural way than the groups $\LB\nn$. For this reason we develop in a parallel fashion: 
\begin{enumerate}[label= \alph*)]
\item the theory of loop braids, for the sake of unifying and completing existent literature;
\item the theory of extended loop braids, because they appear to be the most natural analogue of classical braids, when considering corresponding objects in the $4$-dimensional space.
\end{enumerate}
The paper is organized as follows.  In Section~\ref{S:MCG} we  give a first definition of loop braid groups in terms of mapping classes. 

In Section~\ref{S:configuration} we introduce some configuration spaces of circles in the $3$-ball~$B^3$, and show that their fundamental groups are isomorphic to $\LBE\nn$ and to~$\LB\nn$  (Theorem~\ref{T:PureRing} and Proposition~\ref{P:PureUnrestricted}). This provides a second interpretation of loop braids. We also recall the presentations given in~\cite{BH} for $\LB\nn$ and~$\LBE\nn$, and for the \emph{pure subgroups} of the first family of groups. 

In Section~\ref{S:automorphisms} we introduce subgroups of \emph{permutation-conjugacy automorphisms} $\PC\nn^\ast$ of the group of automorphisms of~ $\F\nn$, the free group of rank~$\nn$, and recall a result of Dahm~\cite{Dahm}, published by Goldsmith~\cite{GOLD}, stating that extended loop braid groups are isomorphic to~$\PC\nn^\ast$~(Theorem~\ref{T:Gold}). This gives us a third interpretation of~$\LBE\nn$ and $\LB\nn$, in terms of automorphisms of~$\F\nn$, and also a characterization of loop braids as automorphisms similar to Artin's characterization of braids. At this point we also give a presentation for the pure subgroups of the groups~$\LBE\nn$ (Proposition~\ref{P:PresNuova}).

Section~\ref{S:ribbon} brings a more topological viewpoint on loop braids, as particular classes of braided surfaces in a $4$-dimensional ball~$B^4$ that we call \emph{(extended) ribbon braids}. We establish an isomorphism between the groups of (extended) ribbon braids and the (extended) loop braid groups (Theorems~\ref{T:flying} and~\ref{T:Extflying}). This is the fourth interpretation of loop braids that we consider.

Sections~\ref{S:welded} is devoted to representations in dimension $3$ and $2$ of loop braids. Here we introduce \emph{broken surfaces} as well as \emph{welded diagrams}. Broken surfaces are projections in general position of loop braids in the $3$-dimensional space. Welded diagrams are $2$-dimensional diagrams, introduced in~\cite{FRR}: passing through broken surfaces and using results of~\cite{BH} we introduce an isomorphism between the groups of these diagrams and the loop braid groups, seen as groups of ribbon braids (Theorem~\ref{T:IsoTube}). 

In Section~\ref{S:Gauss} we discuss a combinatorial description of loop braids through \emph{Gauss diagrams}. These are diagrams, in the spirit of~\cite{PV}, encoding information about a knotted object. This point of view allows us to see loop braid groups as quotients of the groups of virtual braids. This is done using the isomorphisms between loop braid groups and welded braid groups. Though Gauss diagrams has already been used as an equivalent formulation of welded objects (see for example~\cite{Bn, BG, ABMW2}), to formally prove the isomorphism between the groups of welded Gauss diagrams and welded braid diagrams we need to use results from~\cite{Ci} on virtual braids.

Finally Section~\ref{S:other} contains a brief history of these objects, and some other references to find about topological applications and what is known about a representation theory for loop braid groups.

\section{Mapping class groups of a trivial link of unknotted circles in \texorpdfstring{$B^3$}{B\textthreesuperior}}
\label{S:MCG}
In this section we introduce the mapping class group of a $3$-manifold with respect to a submanifold. Then we present a first definition for loop braid groups in terms of mapping classes of a $3$-ball with $n$ circles that are left setwise invariant in its interior. We also introduce three relatives of these groups.
To introduce the main tools we follow \cite{GOLD} and~\cite{KTD}.

Let $M$ be a compact, connected, orientable $3$-manifold, possibly with boundary, and $N$ an orientable submanifold contained in the interior of $M$, not necessarily connected or non-empty. A \emph{self-homeomorphism} of the pair of manifolds $(M,N)$ is an homeomorphism $f \colon M \to M$ that fixes $\partial M$ pointwise, preserves orientation on~$M$, and globally fixes~$N$. Every self-homeomorphism of $(M, N)$ induces a permutation on the connected components of~$N$ in the natural way.

We denote by $\Homeo(M; N)$  the group of self-homeomorphisms of~$(M, N)$ that preserve orientation on both $M$ and~$N$. The multiplication in  $\Homeo(M; N)$ is given by the usual composition. We denote by $\Homeo(M)$ the group $\Homeo(M; \emptyset)$. Moreover we denote by~$\PHomeo(M; N)$ the subgroup of self-homeomorphisms of $(M, N)$ that send each connected component of $N$ to itself. 

We remark that  $\Homeo(M; N)$ is a topological group when equipped with the compact-open topology. The embedding of  $\Homeo(M; N)$ into $\Homeo(M)$ makes  $\Homeo(M; N)$ a closed subgroup of the topological group~$\Homeo(M)$. 

Let $I$ be the unit interval. Two self-homeomorphisms $f_0, f_1$ of $(M, N)$ are said to be \emph{isotopic} if they can be extended to a family $\{f_t\}_{t \in I}$ of self-homeomorphisms of $(M,N)$ such that the map $M \times I \to M$,  sending $(x, t) \to f_t(x)$, is continuous. The isotopy relation is an equivalence relation and isotopic self-isomorphisms induce the same permutation on the connected components of~$N$. 

\begin{defn}
The \emph{mapping class group of a $3$-manifold $M$ with respect to a submanifold~$N$}, denoted by~$\MCG M {N}$, is the group of isotopy classes of self-homeomorphisms of~$\Homeo(M; N)$, with multiplication determined by composition. We denote by $\textrm{MCG}( M)$ the mapping class group~$\MCG M {\emptyset}$.

The  \emph{pure mapping class group of a $3$-manifold $M$ with respect to a submanifold $N$}, denoted by $\PMCG M {N}$, is the subgroup of elements of $\MCG M {N}$ that send each connected component of $N$ to itself.
\end{defn}

\begin{rmk}
\label{R:pi0}
It is known (\cite{Kel}) that a map $f$ from a topological space $X$ to $\Homeo(M; N)$ is continuous if and only if the map $X \times M \to M$ sending $(x, y) \mapsto f(x)(y)$ is continuous. Taking $X$ equal to the unit interval~$I$, we have that two self-homeomorphisms are isotopic if and only if they are connected by a path in~$\Homeo(M; N)$. Therefore~$\MCG M {N} =  \pi_0(\Homeo(M; N))$. In a similar way we have that $\PMCG M {N}= \pi_0(\PHomeo(M; N))$.
\end{rmk}

We consider now a bigger class of self-homomorphisms of the pair~$(M, N)$, removing the condition of preserving orientation on~$N$. We add an ``$^\ast$'' to the notation of the submanifold to indicate this. We denote by $\Homeo(M; N^\ast)$  the group of self-homeomorphisms of~$(M, N)$. The multiplication in  $\Homeo(M; N^\ast)$ is given by the usual composition. Also this group is a closed subgroup of the topological group~$\Homeo(M)$. The isotopy relation considered remains the same. 
We denote by $\PHomeo(M; N^\ast)$ the subgroup of $\Homeo(M; N^\ast)$ of self-homeomorphisms that send each connected component of $N$ to itself.

\begin{defn}
The \emph{extended mapping class group of a $3$-manifold $M$ with respect to a submanifold~$N$}, denoted by~$\MCG M {N^\ast}$, is the group of isotopy classes of self-homeomorphisms of~$\Homeo(M; N^\ast)$, with multiplication determined by composition. 

The  \emph{pure extended mapping class group of a $3$-manifold $M$ with respect to a submanifold $N$}, denoted by $\PMCG M {N^\ast}$, is the subgroup of $\MCG M {N^\ast}$ of elements that send each connected component of $N$ to itself.
\end{defn}

As in Remark \ref{R:pi0}, we have an equivalent definition of the groups $\MCG M {N^\ast}$ and $\PMCG M {N^\ast}$ in terms of connected components of subgroups of~$\Homeo(M)$. We have that $\MCG M N$ can be defined as $\pi_0(\Homeo(M; N^\ast))$, and $\PMCG M {N^\ast}$ as~$\pi_0(\PHomeo(M; N^\ast))$.

\begin{exam}
\label{E:MappingBraid}
Fix $\nn \geq 1$. Let us take $M$ to be the disk $D^2$, and $N$ to be a set of $\nn$ distinct points $P =\{p_1, \dots, p_\nn\}$ in the interior of~$D^2$. In this case we have that  $\PMCG {D^2} {P}$ and $\PMCG{ D^2} {P^\ast}$ are isomorphic, since there is no choice on the orientation of a point. For the same reason   $\MCG{ D^2} {P}$ is isomorphic to $\MCG{ D^2} {P^\ast}$. Then the group $\PMCG{ D^2} P$ is isomorphic to $\PB\nn$, the \emph{pure braid group on $n$ strands}, and  $\MCG{ D^2} {P}$ is isomorphic to $\BB\nn$, the \emph{braid group on $n$ strands}, as defined, for example, in~\cite{GM} or~\cite[Chapter~1.6]{KTD}.
\end{exam}

We conclude this section with the definition of the main objects we are interested in: the loop braid group, the extended loop braid group, and their respective pure subgroups.

\begin{defn}
\label{D:LBN}
Let us fix $\nn \geq 1$, and let $C = C_1 \sqcup \cdots \sqcup C_\nn$ be a collection of $n$ disjoint, unknotted, oriented circles, that form a trivial link of $\nn$ components in~$\Rr^3$. The exact position of $C$ is irrelevant because of~\cite[Corollary~3.8]{GOLD}. So in the following we assume that $C$ is contained in the $xy$-plane in the $3$-ball~$B^3$. The \emph{loop braid group on $n$ components}, denoted  by $\LB\nn$, is the mapping class group~$\MCG{B^3}{C}$. The \emph{pure loop braid group on $n$ components},  denoted by~$\PLB\nn$,  is the pure mapping class group~$\PMCG{B^3}C$. In a similar way the \emph{extended loop braid group}, denoted  by $\LBE\nn$, is the extended mapping class group~$\MCG{B^3}{C^\ast}$. The \emph{pure extended loop braid group},  denoted by~$\PLBE\nn$,  is the pure extended mapping class group~$\PMCG{B^3}{C^\ast}$.
\end{defn}

This definition appears in~\cite{GOLD}, in terms of \emph{motion groups}, which we have reformulated in terms of mapping class groups. We note in particular that the groups $\LBE\nn$ coincide with the motion groups of a trivial unlink defined in~\cite{GOLD}. However the nomenclature ``loop braid groups'' is due to Baez, Crans, and Wise~\cite{BWC}.

\section[The configuration space of \texorpdfstring{$\nn$}{n} circles]{The configuration space of a trivial link of unknotted circles}
\label{S:configuration}
The second interpretation of loop braid groups $\LB\nn$ that we give is in terms of configuration spaces, and has been introduced in~\cite{BH} . We recall some notions and results about configuration space. Then we exhibit the isomorphism from the fundamental groups of certain configuration spaces to the loop braid groups $\LB\nn$. At the end of this section we will give presentations for the groups $\LB\nn$ and~$\LBE\nn$.

\begin{defn}
\label{D:configurations}
Let $\nn \geq 1$, and let $\Ri\nn$ be the space of configurations of $\nn$ Euclidean, unordered, disjoint, unlinked circles in $B^3$. The \emph{ring group} $\R\nn$ is its fundamental group. Let $\UnRi\nn$ be the space of configurations of $\nn$ Euclidean, unordered, disjoint, unlinked circles in $B^3$ lying on planes parallel to a fixed one. The \emph{untwisted ring group} $\UR\nn$ is its fundamental group. 

Similarly, let $\PRi\nn$ be the space of configurations of $\nn$ Euclidean ordered, disjoint, unlinked circles. The \emph{pure ring group} $\PR\nn$ is its fundamental group. Finally, let $\PUnRi\nn$ be the space of configurations of $\nn$  Euclidean, ordered, disjoint,  unlinked circles lying on planes parallel to a fixed one. The \emph{pure untwisted ring group} $\PUR\nn$ is its fundamental group. 
\end{defn}

\begin{rmk}
\label{R:connect}
The path connectedness of $\Ri\nn$ is proved in~\cite{FreedmanSkora}.
\end{rmk}

Let $S_\nn$ be the group of permutations of an $\nn$-elements set, and $p$ the orbit projection $\PRi\nn \to \Ri\nn$ that forgets the order of the circles. We observe that $p$ is a regular $n!$-sheeted cover (see Remark~\ref{R:connect}), with $S_\nn$ as group of deck transformations. This covering is associated to~$p_\ast (\pi_1(\PRi\nn))$. From this follows that $\PR\nn$ is a subgroup of $\R\nn$, and we have the short exact sequence
\[ \begin{CD}
1 @>>> \PR\nn @>>> \R\nn @>>> S_\nn @>>> 1.
\end{CD} \]

Similarly, for $\PUR\nn$ and~$\UR\nn$ we can consider the short exact sequence
\[
\begin{CD}
1 @>>> \PUR\nn @>>> \UR\nn @>>> S_\nn @>>> 1. 
\end{CD}
\]

\begin{rmk}
Let $\mathcal{M}_n$ be the configuration space of $\nn$ ordered distinct points in the complex plane~$\Cc$, \ie, the set of $\nn$-tuples $(z_1, \dots, z_n)$ such that $z_i \neq z_j $ for $i \neq j$. It is easy to show that 
its fundamental group is isomorphic to the pure braid group on $\nn$ strands~$\PB\nn$, 
see for instance~\cite[Chapter~1.4]{KTD}. In the same way, the fundamental group of the configuration space of $\nn$ unordered distinct points in~$\Cc$, meaning the fundamental group of the quotient of $\mathcal{M}_n$ by the action of~$S_\nn$, is one definition of the braid group on $\nn$ strands~$\BB\nn$.
\end{rmk}

The following result of Brendle and Hatcher shows the relation between the untwisted ring groups $\UR\nn$ and the ring groups~$\R\nn$.

\begin{prop}[{\cite[Proposition~2.2]{BH}}]
For $\nn \geq 1$, the natural map $\UR\nn \to \R\nn$ induced by the inclusion $\UnRi\nn \to \Ri\nn$ is injective.
\end{prop}

In order to prove that the groups $\R\nn$  are isomorphic to the loop braid groups~$\LB\nn$, seen as the mapping class groups~$\MCG {B^3} {C^\ast}$, we need some results that we will list below.

\begin{thm}[{\cite[Theorem~1]{BH}}]
\label{T:RelaxingCircles}
For $\nn \geq 1$, the inclusion of $\Ri\nn$ into the space of configurations of all smooth trivial links of $\nn$ components in $\Rr^3$, denoted by~$\mathcal{L}_{\nn}$, is a homotopy equivalence.
\end{thm}
The result allows us to consider the fundamental group of the configuration space of smooth trivial links as isomorphic to~$\R\nn$. 
Let $\mathcal{PL}_{\nn}$ be the space of configurations of all smooth trivial links of $\nn$ ordered components in~$\Rr^3$. From Theorem~\ref{T:RelaxingCircles} we can deduce the following corollary about ordered condifuration spaces.

\begin{cor}
For $\nn \geq 1$, the inclusion of $\PRi\nn$ into the space of configurations of all ordered smooth trivial links of $\nn$ components in $\Rr^3$,  denoted by~$\mathcal{PL}_{\nn}$, is a homotopy equivalence.
\end{cor}
\proof
 We remark that the proof of Theorem~\ref{T:RelaxingCircles} is carried on locally, by considering each component in a sphere, disjoint from the other component's spheres, and rounding each component of the configuration. Then, attaching the order information on the components is left unaffected by the transformation.
\endproof

\begin{rmk}
In Section~\ref{S:MCG} we defined the loop braid group $\LB\nn$ (and the pure, extended, and pure extended versions) as the \emph{topological} mapping class group of a $3$-ball $B^3$ with respect to a collection $C$ of $n$ disjoint, unknotted, oriented circles, that form a trivial link of $\nn$ components. This means that we defined it in terms of self-homeomorphisms. 

However, as Wilson recalls in~\cite{Wilson}, it follows from two results of Wattenberg~\cite{Wattenberg}[Lemma 1.4 and Lemma 2.4] that the topological mapping class group of the $3$-space with respect to a collection of  $n$ disjoint, unknotted, oriented circles, that form a trivial link, is isomorphic to the $C^\infty$ mapping class group, defined in terms of diffeomorphisms.

A \emph{self-diffeomorphism} of the pair of manifolds $(B^3,C)$ is a diffeomorphism $f \colon B^3 \to B^3$ that fixes $\partial B^3$ pointwise, preserves orientation on~$B^3$, and fixes setwise~$C$.

We denote by $\Diffeo(B^3; C)$  the group of self-diffeomorphisms of~$(B^3, C)$ that preserve orientation on both $B^3$ and~$C$. We denote by $\Diffeo(B^3)$ the group $\Diffeo(B^3; \emptyset)$. Moreover we denote by~$\PDiffeo(B^3; C)$ the subgroup of self-diffeomorphisms of $(B^3, C)$ that send each connected component of $C$ to itself. In the same spirit of Section~\ref{S:MCG}, we can define  $\Diffeo(B^3; C^\ast)$ and its subgroup~$\PDiffeo(B^3; C^\ast)$. 

Relying on the results of Wattenberg, we have that
\[ \pi_0(\Homeo(B^3; C)) \cong \pi_0(\Diffeo(B^3; C)) \ \mbox{and} \ \pi_0(\Homeo(B^3; C^\ast)) \cong \pi_0(\Diffeo(B^3; C^\ast)).  \]
\end{rmk}


We take $C = (C_{1}, \dots, C_{n})$ to be an ordered tuple of $\nn$ disjoint, unlinked, trivial circles living in~$B^3$.
We consider the space of configurations of ordered smooth trivial links of $\nn$ components $\mathcal{PL}_\nn$ with the topology induced by the distance defined as follows. Each trivial knot of the configuration can be seen in a $3$-ball $B^3$, and bring smooth, it admits parametrizations.  Let us fix a distance $d$ on~$B^3$. Then taken two trivial knots $C_1$ and $C_2$ in $B^3$, we denote by $d(C_1, C_2)$ the $\min  \{\max \{d(p_1(t), p_2(t))\}\} $, where the $\min$ is considered on the parametrizations $p_i \colon  S^1 \to B^3$ with $i=1, 2$ for $C_1$ and $C_2$, and the $\max$ is considered on the parameter~$t \in S^1$. 

We define an \emph{evaluation map} 
\begin{equation}
\label{E:evaluation}
\ep \colon \Diffeo(B^3) \longrightarrow \mathcal{PL}_{\nn}
\end{equation}
sending a self-diffeomorphism $f$ to~$f(C)$. Remark that $f(C)$ is an ordered tuple of $\nn$ disjoint, unlinked, trivial, smooth knots living in~$B^3$, since the diffeomorphism could have deformed the circles of~$C$. This map is surjective and continuous for construction.

\begin{lem}
\label{L:fibration}
For $\nn \geq 1$, the evaluation map $\ep \colon \Diffeo(B^3) \to \mathcal{PL}_{\nn}$ is a locally trivial fibration, with fibre~$\PDiffeo(B^3; C^\ast)$.
\end{lem}
\proof
Let us consider a point  $C^0 = (C_{1}^0, \dots, C_{n}^0)$ in~$\mathcal{PL}_{\nn}$, \ie, an ordered tuple of $\nn$ disjoint, smooth, unlinked, trivial knots living in~$B^3$. Note that:
\[
 \ep\inv (C^0)=\{f \in \Diffeo (B^3) \mid f(C_{i}^0)=C_{i}^0 \mbox{ for } i=1, \dots, n \} \cong \PDiffeo(B^3; {C^0}^\ast).
\]
Moreover, for any $C$ in $\mathcal{PL}_\nn$ there is an isomorphism between the groups $\PDiffeo(B^3; {C^0}^\ast)$ and~$\PDiffeo(B^3; {C}^\ast)$. This means that $\PDiffeo(B^3; {C}^\ast)$ is the fibre of the fibration.

We have already remarked that $\ep$ is surjective. For $\ep$ to be a locally trivial fibration we need to prove that for every point $C \in \mathcal{PL}_{\nn}$ there is a neighbourhood ${U}_{C} \subset \mathcal{PL}_{\nn}$ together with a homeomorphism from  $U_{C} \times \PDiffeo(B^3; {C}^\ast)$ to $\ep\inv (U_{C}) $ whose composition with $\ep$ is the projection to the first factor~$U_{C} \times \PDiffeo(B^3; {C}^\ast) \to U_{C}$. 
This can be done concretely constructing the local product structure.

However an alternative proof of this Lemma consists in the following remark. The topological group $\Diffeo(B^3)$ acts transitively on the disc in the sense that: if $(C_1, \dots, C_\nn)$ is a collection of $\nn$ circles with the usual conditions, and $(\gamma_1, \dots, \gamma_\nn)$ is another collection with the same conditions, then there is a diffeomorphism $h$ of $\Diffeo(B^3)$ such that $h(C_i)=\gamma_i$ for all $i \in\{1, \dots, \nn\}$. As already remarked, if $h$ is an element of $\Diffeo(B^3, C^\ast)$, then $h(C_i)=C_i$ for all $i \in\{1, \dots, \nn\}$, and if $h, h'$ are such that~$\ep(h)=\ep(h')$, then they are in the same left coset of $\Diffeo(B^3, C^\ast)$ in~$\Diffeo(B^3)$. This means that $\PDiffeo(B^3, C^\ast)$ has a local cross-section in $\Diffeo(B^3)$ with respect to~$\ep$. Then, for~\cite{STE}[Section~7.4], we have the result. 
\endproof


To prove the next theorem, we also need a result from Hatcher on the group~$\Diffeo(B^3)$.

\begin{thm}[{\cite{Hatcher}, Appendix}]
\label{T:Diffeo3ball}
The group $\Diffeo(B^3)$ is contractible. 
\end{thm}

\begin{thm}
\label{T:PureRing}
For $\nn \geq 1$, there are natural isomorphisms between the pure ring group $\PR\nn$ and the pure extended loop braid group~$\PLBE\nn$, and between their respective unordered versions $\R\nn$ and~$\LBE\nn$.
\end{thm}

\proof
Let $\ep$ be the evaluation map \eqref{E:evaluation}. For Lemma~\ref{L:fibration}, we have the short exact sequence
\[
\begin{CD}
1 @>>> \PDiffeo(B^3; C^\ast) @>i>> \Diffeo(B^3) @>{\ep}>>  \mathcal{PL}_{\nn} @>>> 1.
\end{CD}
\]
This induces a long exact sequence of homotopy groups: 
\[
\minCDarrowwidth15pt
\begin{CD} 
\cdots @>>> \pi_1\big(\Diffeo(B^3)\big) @>{\ep_\star}>> \pi_1 \big(\mathcal{PL}_{\nn} \big) @>\partial>> \pi_0 \big(\PDiffeo(B^3; C^\ast)\big) @>{i_\star}>> {} \\
{}@. {}@. {} @>{i_\star}>> \pi_0 \big(\Diffeo(B^3)\big) @>>> \cdots
\end{CD}
\]

Both the groups  $\pi_1\big(\Diffeo(B^3)\big)$ and $\pi_0 \big(\Diffeo(B^3)\big)$ are trivial: both of these statements follow from Theorem~\ref{T:Diffeo3ball}. We recall that  $\pi_1 (\mathcal{PL}_{\nn}) $ is isomorphic to $\PR\nn $ and $\pi_0 \big(\PDiffeo(B^3; C^\ast)\big)$ is isomorphic to $\PLBE\nn$. Then we have an isomorphism between $\PR\nn$ and~$\PLBE\nn$.

Let us now consider $\R\nn$. We can construct the following commutative diagram: 
\[
\begin{CD}
1 @>>> \PLBE\nn @>>> \LBE\nn @>>> S_\nn @>>> 1 \\
@.		@V{\cong}VV		@VVV		 @| 		@.  \\
1 @>>>  \PR\nn @>>> \R\nn @>>> S_\nn @>>> 1 .
\end{CD}
\]
The bijectivity of the central homomorphism follows from the five lemma. 
\endproof

\begin{rmk}
The proof of Theorem~\ref{T:PureRing} can be simplified in a way that does not require Theorem~\ref{T:Diffeo3ball}. Remark that $\ep$ is homotopic to the constant map. In fact the collection of trivial knots $C$ in $B^3$ can be deformed by isotopy to a collection of trivial knots in~$B^3_{x \geq 0}$, the right half of the $3$-ball~$B^3$. Moreover an element $f \in \Diffeo(B^3)$ coincides with the identity on~$\partial B^3$. Thus, it can be deformed to be the identity on~$B^3_{x \geq 0}$. 
The composition of these two deformations is a homotopy $\Diffeo(B^3) \times I \to \mathcal{PL}_{\nn}$ that deforms $\ep$ to the constant map. 

Since the constant map induces the zero map on homotopy groups, whe can consider the following short exact sequence: 
\[
\minCDarrowwidth15pt
\begin{CD} 
1 @>>> \pi_1 \big(\mathcal{PL}_{\nn} \big) @>\partial>> \pi_0 \big(\PDiffeo(B^3; C^\ast)\big) @>{i_\star}>>  \pi_0 \big(\Diffeo(B^3)\big) @>>> 1
\end{CD}
\]
Since $\Diffeo(B^3)$ is path connected~\cite{Cerf}, the group $\pi_0 \big(\Diffeo(B^3)\big)$ is trivial, and we have the isomorphism between $\pi_1 (\mathcal{PL}_{\nn}) $ and~$\pi_0 \big(\PDiffeo(B^3; C^\ast)\big)$.
\end{rmk}

\begin{prop}
\label{P:PureUnrestricted}
For $\nn \geq 1$, there are natural isomorphisms between pure untwisted ring group $\PUR\nn$ and the pure loop braid group~$\PLB\nn$, and between their respective unordered versions $\UR\nn$ and~$\LB\nn$.
\end{prop}
\proof
The group $\PLB\nn$ injects as a normal subgroup in the group~$\PLBE\nn$. In particular, we recall that $\PLB\nn$ is the subgroup of elements of the mapping class group $\PLBE\nn$ that preserve orientation on the $\nn$ connected components of the submanifold $C$ of $B^3$. It is the kernel of the map  $\PLBE\nn \to  \Zz^\nn_2$ sending an homeomorphism reversing the orientation on the $i$th component of $C$, and preserving the orientation on all the other components, to $(0, \dots, 1, \dots, 0)$, where the non-zero entry is in position~$i$. We have the following short exact sequence.
\begin{equation}
 \raisebox{-0.3\height}{\includegraphics{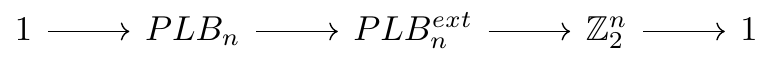}}
\end{equation}
On the other hand, from \cite[Proposition~2.2]{BH}, we have the short exact sequence
\begin{equation}
\label{E:sequenzaSplitting}
 \raisebox{-0.3\height}{\includegraphics{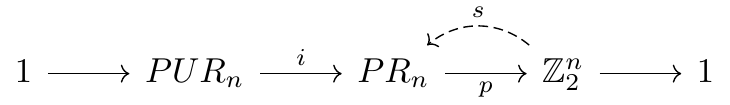}}
\end{equation}
which has an obvious splitting~$s$, obtained by rotating the circles within disjoint balls. From the isomorphism between $\PR\nn$ and $\PLBE\nn$ (Theorem~\ref{T:PureRing}) we get the diagram: 
\[
\begin{CD}
1 @>>>  \PLB\nn @>>> \PLBE\nn @>>> \Zz^n_2  @>>> 1 \\
@.		  @VVV	@V{\cong}VV         @| 	@.  \\
1 @>>>  \PUR\nn @>>> \PR\nn @>>> \Zz^n_2 @>>> 1 .
\end{CD}
\]
For the universal property of the kernel, the first vertical arrow is an isomorphism, which proves the first part of the statement.
To prove the second isomorphism, we consider the commutative diagram
\[
\begin{CD}
1 @>>>  \PLB\nn @>>> \LB\nn @>>> S_\nn @>>> 1 \\
@.		@V{\cong}VV		@VVV		 @| 		@.  \\
1 @>>>  \PUR\nn @>>> \UR\nn @>>> S_\nn @>>> 1 .
\end{CD}
\] 
The bijectivity of the central homomorphism follows from the five lemma. 
\endproof

\begin{rmk}
Note that to prove the isomorphism of Proposition~\ref{P:PureUnrestricted} for loop braids we need to prove the result on extended loop braids first. This is because the first one relies on the fibration on $\mathcal{PL}_{\nn}$ of Lemma~\ref{L:fibration}. Indeed, in the context of configuration spaces, we cannot define a notion of ``preserving orientation on the circles'' when working with a single configuration of $\nn$ circles as a point in the topological space of configurations.
\end{rmk}

Brendle and Hatcher, in \cite[Proposition~3.3]{BH}, also give a presentation for the untwisted ring groups~$\UR\nn$, which are isomorphic to the loop braid groups~$\LB\nn$.

\begin{prop}
\label{P:BrendleHatcher}
For $\nn \geq 1$, the group $\UR\nn$ admits a presentation given by generators $\{\sig\ii, \rr\ii \mid \ii=1, \dots , \nno \}$, subject to relations:
\begin{equation}
\label{E:Upresentation}
\begin{cases}
\sig{i} \sig j = \sig j \sig{i}  \, &\text{for } \vert  i-j\vert > 1\\
\sig{i} \sig {i+1} \sig{i} = \sig{i+1} \sig{i} \sig{i+1} \, &\text{for } i=1, \dots, \nn-2 \\
\rr{i} \rr j = \rr j \rr{i}  \, &\text{for }  \vert  i-j\vert > 1\\
\rr\ii\rr{i+1}\rr\ii = \rr{i+1}\rr\ii\rr{i+1} \, &\text{for }  i=1, \dots, \nn-2  \\
\rrq{i}2 =1 \, &\text{for }  i=1, \dots, \nno \\
\rr{i} \sig{j} = \sig{j} \rr{i}   \, &\text{for }  \vert  i-j\vert > 1\\
\rr{i+1} \rr{i} \sig{i+1} = \sig{i} \rr{i+1} \rr{i} \,  &\text{for }  i=1, \dots, \nn-2  \\
\sig{i+1} \sig{i} \rr{i+1} = \rr{i} \sig{i+1} \sig{i}  \,  &\text{for }  i=1, \dots, \nn-2. \\
\end{cases}
\end{equation}
\end{prop}

From this proposition it follows:
\begin{cor}
\label{C:PresLB}
For $\nn \geq 1$, the group $\LB\nn$ admits the presentation given in Proposition~\ref{P:BrendleHatcher}.
\end{cor}

In \cite[Proposition~3.7]{BH} they also obtain a presentation for $\R\nn$, adding to the presentation in Proposition~\ref{P:BrendleHatcher} generators $\{\tau_\ii \mid \ii=1, \dots , \nn \}$, and some relations.

\begin{prop}
\label{P:PresRn}
For $\nn \geq 1$, the group $\R\nn$ admits a presentation given by generators $\{\sig\ii, \rr\ii \mid \ii=1, \dots , \nno \}$ and $\{\tau_\ii \mid \ii=1, \dots , \nn \}$, subject to relations:
\begin{equation}
\label{E:Rpresentation}
\begin{cases}
\tau_{i} \tau_j = \tau_j \tau_{i}  \, &\text{for }    i \neq j \\
\tau_\ii^2=1 \, &\text{for }  i=1, \dots, \nn \\
\sig{i} \tau_j = \tau_j \sig{i}  \, &\text{for }  \vert  i-j\vert > 1 \\
\rr{i} \tau_j = \tau_j \rr{i}  \, &\text{for }  \vert  i-j\vert > 1 \\
\tau_i \rr\ii = \rr\ii \tau_{i+1} \, &\text{for }  i=1, \dots, \nno  \\
\tau_\ii \sig\ii = \sig\ii \tau_{i+1}  \, &\text{for } i=1, \dots, \nno \\
\tau_{i+1} \sig\ii = \rr\ii \siginv\ii \rr\ii \tau_\ii  \, &\text{for } i=1, \dots, \nno. \\
\end{cases}
\end{equation}
\end{prop}

From this proposition it follows:
\begin{cor}
\label{C:PresLBE}
For $\nn \geq 1$, the group $\LBE\nn$ admits the presentation given in Proposition~\ref{P:PresRn}.
\end{cor}

\begin{figure}[hbtp]
\centering
\includegraphics[scale=.6]{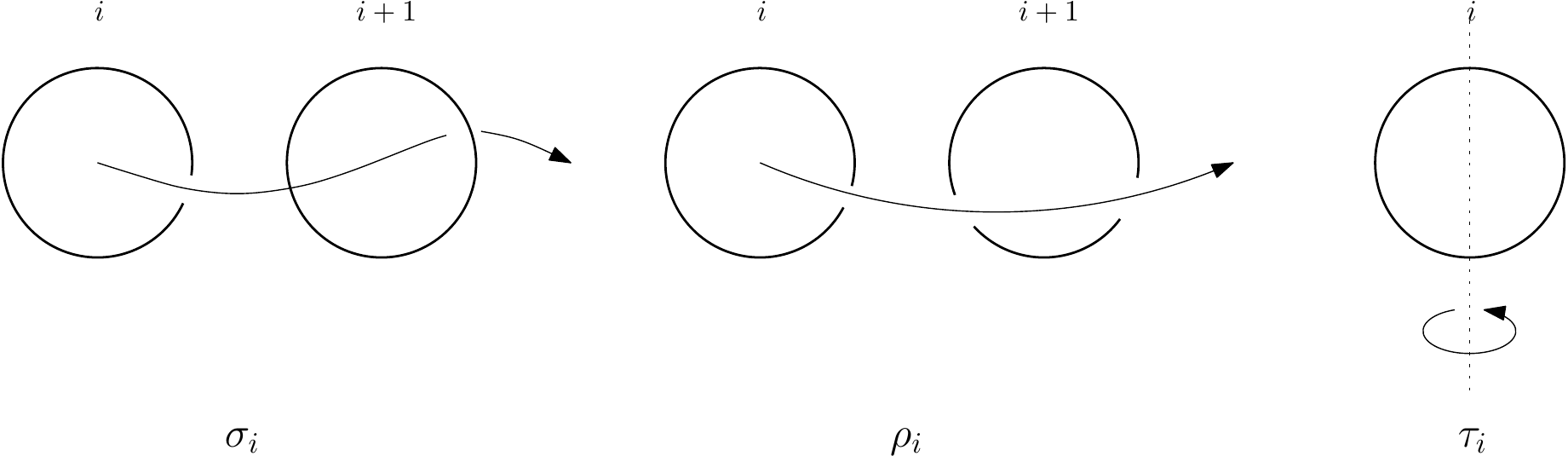}
\caption{Elements $\sig\ii$, $\rr\ii$ and $\tau_\ii$.}
\label{F:Flips}
\end{figure}

The elements~$\sig\ii$,~$\rr\ii$, and $\tau_\ii$ of the presentation represent the following loops in $\R\nn$: if we place the $\nn$ rings in a standard position in the $yz$-plane with centers along the $y$-axis, then the $\sig\ii$ is the loop that permutes the $\ii$-th and the $(\ii+1)$-st circles by passing the $\ii$-th circle through the $(\ii+1)$-st; $\rr\ii$ permutes them passing the $\ii$-th around the $(\ii+1)$-st, and $\tau_i$ changes the orientation (``flips'') the $i$-th circle, see Figure~\ref{F:Flips}\footnote{We reverse the notations used in \cite{BH} for $\sig\ii$ and~$\rr\ii$, see also, for example, notations used in \cite{FRR, GOLD, BBVW}}.

\begin{figure}[hbtp]
\centering
\includegraphics[scale=.6]{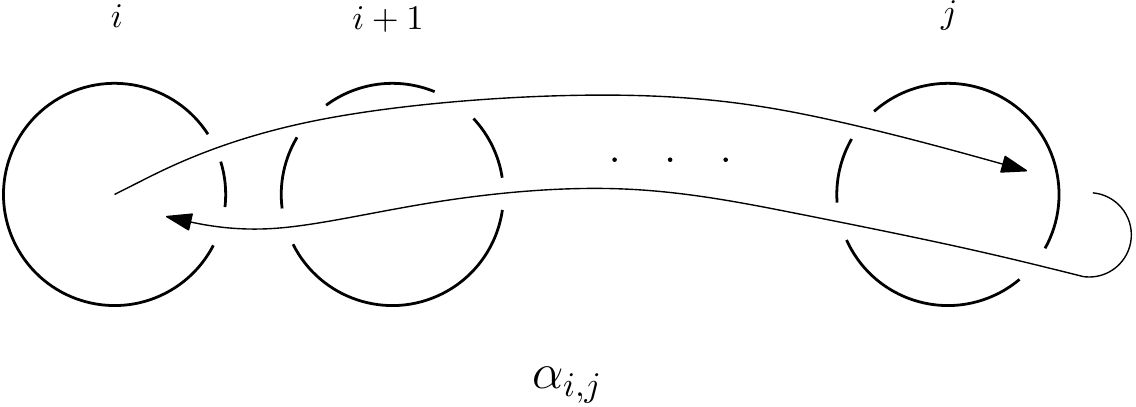}
\caption{Element $\alpha_{ij}$.}
\label{F:Alpha}
\end{figure}

Let $\alpha_{ij}$ be the elements of~$\PUR\nn$, representing the movement of the $i$-th circle passing throug the $j$-th circle and going back to its position (see Figure~\ref{F:Alpha}). Brendle and Hatcher also give a presentation for~$\PUR\nn$, the pure subgroups of the groups~$\UR\nn$. 

\begin{prop}
\label{P:BHpresentazionepura}
For $\nn \geq 1$, the group $\PUR\nn$ admits a presentation with generators $\alpha_{ij}$ for $1\leq i \neq j \leq \nn$ and relations:
\begin{equation}
\begin{cases}
\alpha_{ij} \alpha_{kl}= \alpha_{kl} \alpha_{ij} \\
\alpha_{ik} \alpha_{jk}=\alpha_{jk} \alpha_{ik} \\
\alpha_{ij} (\alpha_{ik} \alpha_{jk}) = (\alpha_{ik} \alpha_{jk}) \alpha_{ij}.
\end{cases}
\end{equation}
\end{prop}

\begin{cor}
\label{C:PresPLB}
For $\nn \geq 1$, the group $\PLB\nn$ admits the presentation given in Proposition~\ref{P:BHpresentazionepura}.
\end{cor}

In Corollaries~\ref{C:PresLB},~\ref{C:PresLBE}, and~\ref{C:PresPLB} we exhibited presentations for the loop braid groups~$\LB\nn$, for the extended loop braid groups~$\LBE\nn$, and for the pure loop braid groups~$\PLB\nn$. At the end of next section we will also give a presentation for the pure extended loop braid groups~$\PLBE\nn$.

\section{Conjugating automorphisms of the free group}
\label{S:automorphisms}
In this section we give an interpretation of $\LB\nn$ and $\LBE\nn$ in terms of automorphisms of~$\F\nn$, the free groups of rank~$\nn$. 
Fixing $\nn \geq 1$, we consider automorphisms that send each generator of $\F\nn$ to a conjugate of some generator: these are in bijection with elements of~$\LB\nn$. When considering automorphisms that send each generator to a conjugate of some generator, or of the inverse of some generator, we have a bijection with elements of $\LBE\nn$. Finally, when considering automorphisms that send a generator of $\F\nn$ to a conjugate of itself, we have a bijection with elements of the pure loop braid group~$\PLB\nn$, and if a generator is sent to the conjugate of itself or its inverse, then we have a bijection with elements of~$\PLBE\nn$.
We start from a result of Dahm's unpublished thesis~\cite{Dahm}, that appears in the last section of Goldsmith's paper~\cite{GOLD}. We remark that the result is established for~$\LBE\nn$, and that the result for $\LB\nn$ is a consequence of it.

\begin{thm}[{\cite[Theorem~5.3]{GOLD}}]
\label{T:Gold}
For $\nn \geq 1$, there is an injective map from the extended loop braid group $\LBE\nn$ into~$\Aut(\F\nn)$, where $\F\nn$ is the free group on $\nn$ generators $\{x_1, \dots, x_\nn\}$, and its image is the subgroup~$\PC\nn^\ast$, consisting of all automorphisms of the form $\alpha \colon x_\ii \mapsto w_\ii\inv x^{\pm 1}_{\pi({\ii})} w_\ii$ where $\pi$ is a permutation and $w_\ii$ is a word in~$\F\nn$. Moreover, the group $\PC\nn^\ast$ is generated by the automorphisms $\{ \sig1, \dots \sig\nno,\rr1, \dots \rr\nno, \tau_1, \dots, \tau_\nn\}$ defined as: 
\begin{align}
\label{E:sigma}
\sig\ii &: \begin{cases}
            x_\ii \mapsto x_{\ii+1}; &\\
            x_{\ii+1} \mapsto x_{\ii+1}\inv x_\ii x_{\ii+1}; &\\
            x_\jj \mapsto x_\jj, \ &\text{for} \ \jj \neq \ii, \ii+1.
        \end{cases} \\
\label{E:rho}
\rr\ii &: \begin{cases}
            x_\ii \mapsto x_{\ii+1}; & \\
            x_{\ii+1} \mapsto x_\ii; &\\
            x_\jj \mapsto x_\jj, \ &\text{for} \ \jj \neq \ii, \ii+1.
		\end{cases} \\
\label{E:tau}
\tau_\ii &: \begin{cases}
            x_\ii \mapsto x\inv_\ii; &\\
            x_\jj \mapsto x_\jj, \ &\text{for} \ \jj \neq \ii.
        \end{cases}
\end{align}
\end{thm}


\begin{rmk}
\label{R:AbuseNotation}
The correspondence between the elements of the mapping class group representing the movements $\sig\ii$, $\rr\ii$ and $\tau_i$ in $\R\nn$, and the automorphisms of $\PC\nn^\ast$ with the same names justifies the abuse of notation. 
\end{rmk}

The following result, which is a consequence of Theorem~\ref{T:Gold}, establishes an isomorphism between $\LB\nn$ and $\PC\nn$, the groups of permutation-conjugacy automorphisms. These groups consist of all automorphisms of the form $\alpha \colon x_\ii \mapsto w_\ii\inv x_{\pi({\ii})} w_\ii$ where $\pi$ is a permutation and $w_\ii$ is a word in~$\F\nn$.

\begin{cor}
\label{C:LoopAutomorphisms}
For $\nn \geq 1$, there is an injection from $\LB\nn$ to~$\Aut(\F\nn)$, where $\F\nn$ is the free group on $\nn$ generators $\{x_1, \dots, x_\nn\}$, and its image is the subgroup~$\PC\nn$. Moreover, the group $\PC\nn$ is generated by automorphisms~\eqref{E:sigma} and~\eqref{E:rho}.
\end{cor}

\begin{rmk}
\label{R:AutoClassicBraids}
The elements $\sig\ii$ in $\PC\nn$ generate the \emph{braid subgroup} $\BB\nn$ of $\Aut(\F\nn)$ which is well known to be isomorphic to the classical braid group on $\nn$ strings, and the elements $\rr\ii$ generate the \emph{permutation subgroup} $\bar{S}_\nn$ of $\Aut(\F\nn)$ which is a copy of the symmetric group~$S_\nn$.  Moreover Artin provided (see for instance~\cite[Theorem 5.1]{LH}) a characterization of usual braids as automorphisms of free groups of which Theorem~\ref{T:Gold} is the analogue. Let us recall it: an automorphism $\beta \in \Aut(\F\nn)$ lies in $\BB\nn$ if and only if $\beta $ satisfies the following conditions:
\begin{enumerate}[label=\roman*)]
\item $\beta(x_i) = a_i \, x_{\pi(i)} \, a_i\inv,~~1\leq i\leq n$\; ;
\item $\beta(x_1x_2 \ldots x_n)=x_1x_2 \ldots x_n$\; ,
\end{enumerate}
where $\pi \in S_\nn$
and $a_i \in \F\nn$.
\end{rmk}

\begin{rmk}
\label{R:Braid-Permutation}
Fenn, Rim\'{a}nyi and Rourke, in~\cite{FRR}, consider the subgroups $\BP\nn$ of $\Aut(\F\nn)$ generated by both sets of elements $\{\sig\ii \mid i=1, \dots \nn-1\}$ and $\{\rr\ii \mid i=1, \dots \nn-1\}$. They call these groups the \emph{braid-permutation groups}, and they prove independently from Dahm and Goldsmith that they are isomorphic to the permutation-conjugacy groups~$\PC\nn$, and that they admit the presentation given in Proposition~\ref{P:BrendleHatcher}.
\end{rmk}

\begin{rmk}
The presentation given in Proposition~\ref{P:BHpresentazionepura} coincides with the presentation given by McCool in~\cite{McC} for the groups of \emph{basis-conjugating automorphisms} of the free group. Fixed an $\nn \geq 1$, this is the group generated by the automorphisms:
\begin{equation}
\label{E:epsilon}
\alpha_{\ii \jj} : \left\{
        \begin{array}{lll}
            x_{\ii} \mapsto x_{\jj}\inv x_\ii x_{\jj};\\
            x_k \mapsto x_k, \ \text{for} \ k \neq \ii, \jj.
        \end{array}
    \right.
 \end{equation}
\end{rmk}

\begin{rmk}  
In~\cite{SAV} Savushkina proves that the centers of the groups $\PLB\nn$ are trivial, and gives a presentation for the groups of \emph{permutation-conjugacy automorphisms}, isomorphic to the loop braid groups~$\LB\nn$.
\end{rmk}

To complete the picture we use the conjugating automorphisms point of view to give a presentation for~$\PLBE\nn$, the pure subgroups of the groups of extended loop braids~$\LBE\nn$. 

\begin{prop}
\label{P:PresNuova}
For $\nn \geq 1$, the group $\PLBE\nn$ admits the following presentation:
\begin{equation}
\label{E:PRnPres}
\Big\langle \{  \alpha_{ij} \mid  1 \leq \ii \neq \jj \leq \nn \}    \cup     \{  \tau_\ii \mid  \ii = 1, \dots \nn \}\    \ \big\vert \    R^\prime \Big\rangle
\end{equation}
where $R^\prime$ is the set of relations:
\begin{equation}
\begin{cases}
\alpha_{ij} \alpha_{kl} = \alpha_{kl} \alpha_{ij} \\
\alpha_{ij} \alpha_{kj} = \alpha_{kj} \alpha_{ij} \\
(\alpha_{ij} \alpha_{kj}) \alpha_{ik} = \alpha_{ik} (\alpha_{ij} \alpha_{kj})\\
\tau_\ii^2 = 1 \ \\
\tau_\ii \alpha_{\ii\jj}  = \alpha_{\ii\jj} \tau_\ii \ \\
\tau_\ii \alpha_{\jj k} = \alpha_{\jj k } \tau_\ii \ \\ 
\tau_\ii \alpha_{\jj\ii} \tau_\ii  = \alpha_{\jj\ii} \inv \ 
\end{cases}
\end{equation}
where different letters stand for different indices, and $\alpha_{ij}$ and $\tau_{i}$ correspond to the automorphisms~\eqref{E:epsilon} and~\eqref{E:tau}.
\end{prop}

\proof
By the split short exact sequence~\eqref{E:sequenzaSplitting}, and the isomorphism between $\PR\nn$ and~$\PLBE\nn$, we have that $\PLBE\nn$ is isomorphic to the semidirect product~$\PLB\nn \rtimes_\varphi \Zz_2^\nn$, where $\varphi \colon \Zz_2^\nn \to \Aut(\PLB\nn)$ is the map defined by:
\[\varphi(x) = \varphi_x \colon g \longrightarrow \ii\inv \big( s(x) i(g) s(x\inv) \big).
\]
Considering the presentation $\Gr{ \{ \tau_\ii \mid  \ii = 1, \dots \nn \} \ \big\vert \ 
\tau_\ii^2 = 1}$ for $\Zz_2^\nn$ we have that $\PLBE\nn$ admits the following presentation:
\begin{equation}
\Big\langle (\alpha_{ij})_{ 1 \leq i \neq  j \leq \nn \}}) \cup  (\tau_{i})_{ 1 \leq \ii \leq \nn \}}) \  \big\vert \ R , S, \{ \tau_\ii  \alpha_{jl}  \tau_\ii\inv = \varphi_{\tau_\ii}( \alpha_{jl} ) \} \Big\rangle
\end{equation}
where the $\alpha_{ij}$ are generators of~$\PLB\nn$, the $\tau_i$ are generators of~$\Zz_2^\nn$, and $R$ and $S$ are the respective sets of relations.

The groups $\PLB\nn$ and $\Zz_2^\nn$ can be seen as subgroups of~$\Aut(F_\nn)$. Then $\varphi$ is the action of $\Zz_2^\nn$ on~$\PLB\nn$, so to understand the mixed relation between the $\alpha_{\ii\jj}$ and the $\tau_ k$ it is sufficient to see how $\tau_ k$ act on~$\alpha_{\ii\jj}$, for  $1 \leq \ii \neq \jj \leq \nn$ and~$ k = 1, \dots \nn$.

Both families of automorphisms are defined on generators $x_\ii$ of~$F_\nn$. In particular, $x_l$ with~$l \neq \ii, \jj, k$, is left unvaried by both automorphisms. Hence it is sufficient to consider the cases where $l$ is equal to $\ii$ or $\jj$, and/or $ k$. In other words,  we are only to consider the mutual positions of three indices. Doing explicit calculation for all the possible cases one finds the relations of the statement (these calculations can be found in~\cite{Dam}).
\endproof

\section{Ribbon braids and flying rings}
\label{S:ribbon}
We now introduce \emph{ribbon braids} and \emph{extended ribbon braids}. Using the results from Section~\ref{S:configuration}, we establish isomorphisms between the groups of (extended) ribbon braids and the (extended) loop braid groups. These isomorphisms give us another interpretation of loop braids as topological knotted objects. We also prove that every isotopy of a ribbon braid in $B^3 \times I$ extends to an isotopy of $B^3 \times I$ itself constant on the boundary. 

Let $\nn \geq 1$. We recall and adapt notations and definitions from \cite{AUD} and~\cite{ABMW}. Let $D_1, \dots, D_\nn$ be a collection of disks in the $2$-ball~$B^2$. Let $C_i =\partial D_i$ be the oriented boundary of~$D_i$.
Let us consider the $4$-ball $B^4 \cong B^3 \times I$, where $I$ is the unit interval. For any submanifold $X \subset B^m \cong B^{m-1} \times I$, with~$m=3,4$, we use the following dictionary. To keep the notation readable, here we denote the interior of a topological space by "$\mathrm{int}(\phantom{x} )$", whereas anywhere else it is denoted by "$\mathring{\phantom{x}}$".
\begin{itemize}
\item $\partial_\ep X = X \cap (B^{m-1} \times \{\ep\})$, with~$\ep \in \{0, 1\}$;
\item $\partial_\ast X = \partial X \setminus \Big( \mathrm{int}(\partial_0 X) \sqcup \mathrm{int}(\partial_1 X) \Big)$;
\item $\overset{*}{X}=X \setminus \partial_\ast X$.
\end{itemize}

The image of an immersion $Y \subset X$ is said to be \emph{locally flat} if and only if it is locally homeomorphic to a linear subspace $\Rr^k$ in $\Rr^m$ for some $k \leq m$, except on $\partial X$ and/or $\partial Y$, where one of the $\Rr$ summands should be replaced by $\Rr_+$. Let $Y_1, Y_2$ be two submanifolds of $B^m$. The intersection $Y_1 \cap Y_2 \subset X$ is called \emph{flatly transverse} if and only if it is locally homeomorphic to the transverse intersection of two linear subspaces $\Rr^{k_1}$ and $\Rr^{k_2}$ in $\Rr^m$ for some positive integers $k_1, k_2 \leq m$ except on $\partial X$, $\partial Y_1$ and/or $\partial Y_2$, where one of the $\Rr$ summands should be replaced by $\Rr_+$. In the next definition we introduce the kind of singularities we consider.

\begin{defn}
\label{D:RibbonDisk}
\emph{Ribbon disks} are intersections $D = Y_1 \cap Y_2 \subset \Rr^4$ that are isomorphic to the $2$-dimensional disk,  such that $D \subset \mathring{Y_1}$, $\mathring{D} \subset\mathring{Y_2}$ and $\partial D $ is an essential curve in~$\partial Y_2$.
\end{defn}
These singularities are the $4$-dimensional analogues of the classical notion of ribbon singularities introduces by Fox in~\cite{F}.

\begin{defn}
\label{D:ribbonBraid}
Let $A_1, \dots, A_\nn$ be locally flat embeddings in~$\stackrel{*}{B^4}$ of $\nn$ disjoint copies of the oriented annulus $S^1 \times I$. We say that   
\[b= \bigsqcup_{i\in \{1, \dots, n \}} A_i \]
is a \emph{geometric ribbon braid} if:
\begin{enumerate}
\item \label{bordo} the boundary of each annulus $\partial A_i$ is a disjoint union  $C_i \sqcup C_j$, for $C_i \in \partial_0 B^4$ and for some $C_j \in \partial_1 B^4$. The orientation induced by $A_i$ on $\partial A_i$ coincides with the one of the two boundary circles $C_i$ and~$C_j$;
\item the annuli $A_i$ are fillable, in the sense that they bound immersed $3$-balls $\subset \Rr^4$ whose singular points consist in a finite number of ribbon disks;
\item \label{trans} it is transverse to the lamination $\bigcup_{t \in I} B^3 \times \{t\}$ of~$B^4$, that is: at each parameter~$t$, the intersection between $b$ and $B^3 \times{t}$ is a collection of exactly $\nn$ circles;
\item \label{wen} the orientations of the circles are concordant, at each parameter~$t$, to the orientations of the circles that compose the boundary of the annulus.
\end{enumerate}

The group of \emph{ribbon braids}, denoted by $\rB\nn$ is the group of equivalence classes of geometric ribbon braids up to continuous deformations through the class of geometric ribbon braids fixing the boundary circles, equipped with the natural product given by stacking and reparametrizing. The unit element for this product is the \emph{trivial ribbon braid}~$U= \bigsqcup_{i\in \{1, \dots, n \}} C_i \times [0, 1]$.
\end{defn}

\begin{rmk}
\label{R:wen}
Let us consider  $I$ in $B^4 = B^3 \times I$ as a time parameter. If one of the $n$ circles that we have at each time $t$ makes a half-turn, we have what is called a \emph{wen} on the corresponding component. One can think of a wen as an embedding in $4$-space of a Klein bottle cut along a meridional circle. The last condition of the definition makes sure that there are no wens on the components of a geometric ribbon braid. A detailed treatment of wens can be found in~\cite{KS}.
\end{rmk}

\begin{rmk}
We recall that also the braid group $\BB\nn$ can be defined in an analogous way, as equivalence classes of geometric braids, see~\cite[Chapter~1.2]{KTD}.
\end{rmk}

The following theorem shows that when two ribbon braids are equivalent in the sense of Definition~\ref{D:ribbonBraid}, there is an ambient isotopy of $\Rr^4$ bringing one to the other.
\begin{thm}
\label{T:isotopy}
Every relative isotopy of a geometric ribbon braid in $B^3 \times I$ extends to an isotopy of $B^3 \times I$ in itself constant on the boundary.
\end{thm}
\proof
We follow step by step the proof given for the case of usual braids in~\cite{KTD}. 
Let $b \subset B^4 \cong B^3 \times I$ be a geometric ribbon braid with $\nn$ components, and let us call $T$ the product $B^3 \times I$. Let
\[
F \colon b \times I \longrightarrow T
\]
be an isotopy of~$b$. Thus, for each $s \in I$, the map 
\begin{align*}
F_s \colon b &\longrightarrow T \\
x &\mapsto F(x, s)
\end{align*}
is an embedding whose image is a geometric ribbon braid, and~$F_0 = \id_b$.

We want to define another continuous map $G \colon T \times  I  \to T $ such that for each $s \in I$ the map $G_s \colon T \to T$ sending $x \in T$ to $G(x,s)$ is a homeomorphism fixing $\partial T$ pointwise and extending~$F_s$, and such that $G_0=\id_{T}$.

The first step to construct it is to consider a set $C=\{C_1, \dots , C_\nn\}$ of $\nn$ disjoint, oriented, unlinked, unknotted circles in the interior of $B^3$, and another $3$-ball $D^3$ such that: $C$ is contained in $\mathring{D^3}$, and $\mathring{D^3}$ is contained in $\mathring{B^3}$. Moreover, let $F(b \times I)$ be contained into~$\mathring{D^3} \times I$. Remark that~$C \times \{0, 1\}$ is the disjoint union of the circles that compose the boundary components of~$b$ and of $F_s(b)$ for all $s \in I$.

For any $(s, t) \in I^2$, denote by $f(s, t)$ the unique subset of $\nn$ circles of $\mathring{D^3}$ such that 
\[
F_s(b) \cap (D^3 \times \{t\}) = f(s, t) \times \{t\}.
\] In other words, $f(s, t)$ is composed by the $\nn$ circles of the ribbon braid $F_s(b)$ at the coordinate $\{t\}$ of the foliation $D^3 \times I = \bigcup_{i\in I} D^3 \times \{t\}$.

So we have a continuous map 
\[
f \colon I^2 \longrightarrow \mathcal{UL}_{\nn} \subset \mathcal{L}_{\nn}
\]
where $ \mathcal{UL}_{\nn}$ is the configuration space of $\nn$ disjoint, oriented, unlinked, unknotted circles that lie on parallel planes.

By definition $f(s, 0)=f(s, 1) =C$ for all $s \in I$ and $b = \bigcup_{t \in I} f(0, t) \times \{t\}$. Let us consider now the loop
\begin{align*}
f_0 \colon I &\longrightarrow D^3\\
 t &\mapsto f(0, t)
\end{align*}
that sends the parameter $t$ to the circles composing $b$ at the coordinate~$t$.

Consider the evaluation fibration $\ep \colon \Diffeo(D^3) \to \mathcal{PL}_{\nn}$ from Lemma~\ref{L:fibration}. Composing $\ep$ with the covering map $\mathcal{PL}_{\nn} \to \mathcal{L}_{\nn}$, seeing $\mathcal{L}_{\nn}$ as the orbit space with of the action of the symmetric group of $\mathcal{PL}_{\nn}$, we define a locally trivial fibration
\[
\tilde{\ep}\colon \Diffeo(D^3) \longrightarrow \mathcal{L}_{\nn}.
\]
By the homotopy lifting property of $\tilde{\ep}$, the loop $f_0$ lifts to a path $\hat{f_0}$ as in the following diagram.
\begin{equation*}
 \raisebox{-0.3\height}{\includegraphics{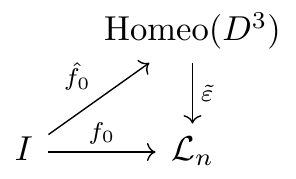}}
\end{equation*}
The lifted map is such that $\hat{f_0}(1)= \id_{D^3}$ and $\hat{f_0}(0)$ is an element of $\Diffeo(D^3; C^\ast)$. Remark that $\tilde{\ep} \circ \hat{f_0}(0)=C=f_0(0)$, $\tilde{\ep} \circ \hat{f_0}(1)=C=f_0(1)$. By the homotopy lifting property of $\ep$, considering the topological pair $(I, \partial I)$, the latter path extends to a lift
\[
\hat{f} \colon I \times I \longrightarrow \Diffeo(D^3)
\]
such that $\hat{f}(s, 1)=\id_{D^3}$ and $\hat{f}(s, 0)=\hat{f}(0, 0)$ for all $s \in I$. 
We define a homeomorphism 
\[
g(s,t) \colon B^3 \longrightarrow B^3
\]
defined by
\begin{equation*}
 g(s,t) = \begin{cases} \id_{B^3} \ &\mbox{on }  B^3 \setminus D^3, \\
						 \hat{f}(s, t) \circ \big(\hat{f}(0, t)\big)\inv  &\mbox{on } x \in D^3.
		   \end{cases}
\end{equation*}

This is a continuous function on $(s, t) \in I^2$ and
\[
g(0, t)=g(s,0)=g(s,1)=\id_{B^3}.
\]
Moreover 
\[
g(s, t)\big(f(0,t)\big)=g(s, t) \big(\hat{f}(0, t)(C))=\hat{f}(s, t)(C)=f(s, t).
\]
Let us define 
\begin{align*}
G\colon T\times I &\longrightarrow T \\
(x, t, s) &\mapsto \big(g(s,t)(x), t\big)
\end{align*}
for $x \in B^3$, $(t, s) \in I \times I$. For each $s \in I$ we have a homeomorphisms
\[
G_s \colon T \longrightarrow T
\]
that fixes $\partial T$ pointwise, extends~$F_s$, and~$G_0=\id_{B^3}$.
\endproof

\begin{rmk}
\label{R:extisotopy}
This result is true also for \emph{surface links}, which are closed surfaces locally flatly embedded in~$\Rr^4$. This is proved with different methods in~\cite[Theorem~6.7]{Kam2}.
\end{rmk} 

\begin{defn}
\label{D:pureRibbon}
A \emph{pure geometric ribbon braid} is a geometric braid as in Definition~\ref{D:ribbonBraid}, for which condition~\eqref{bordo} is replaced with
\begin{enumerate}
\item[\eqref{bordo}$^\prime$] \label{bordoPuro} $\partial A_i = C_i \times \{0, 1\} $ for all $i \in \{1, \dots, n\}$ and the orientation induced by $A_i$ on $\partial A_i$ coincides with that of~$C_i$.
\end{enumerate}
The group of \emph{pure ribbon braids}, denoted by $PrB_n$, is the group of equivalence classes of pure geometric ribbon braids up to continuous deformations through the class of geometric ribbon braids fixing the boundary circles, equipped with the natural product given by stacking and reparametrizing. \end{defn}

\begin{rmk}
The group $PrB_n$ coincides with the kernel of the homomorphism from the group $\rB\nn$ to the group of permutation $S_n$ that associates to a ribbon braid the permutation induced on the boundary circles.
\end{rmk}

\begin{rmk}
\label{R:flying}
Taken $b$ a geometric ribbon braid, the transversality forces $b \cap (B^3 \times \{t\})$ to be the disjoint union of $\nn$ circles, for all $t \in I$. We can though think to a ribbon braid as a trajectory $\beta=\big( C_1(t), \dots, C_\nn(t)\big)$ of circles in~$B^3 \times I$. This interpretation corresponds, in the classical case, to the interpretation of braids as particle dance. A visual idea of this can be found in Dalvit's website~\cite{DAL}. This formulation gives rise to the name \emph{groups of flying rings}, denoted by $FL_n$ and used in~\cite{Bn}.
\end{rmk}

We prove now that, for each $\nn \geq 1$, $PrB_\nn$ and $\rB\nn$ are respectively isomorphic to $\PLB\nn$ and~$\LB\nn$.

\begin{prop}
\label{P:pure}
For $\nn \geq 1$, there is an isomorphism between the pure ribbon braid group $PrB_\nn$ and the pure loop braid group~$\PLB\nn$.
\end{prop}
\proof
Let $\beta=\big( C_1(t), \dots, C_\nn(t)\big)$ be an element of~$PrB_\nn$ described by a parametrization as in Remark~\ref{R:flying}, and set
\[
\varphi(\beta) \colon [0, 1] \longrightarrow \PLB\nn
\]
as the morphism defined by $ t \mapsto \big( C_1(t), \dots, C_\nn(t)\big)$. By definition, $\varphi(\beta)$ is a loop in the configuration space $\PUR\nn$ (Definition~\ref{D:configurations}), and corresponds to an element of $\PLB\nn$ through the isomorphism in Proposition~\ref{P:PureUnrestricted}. This map induces a bijection 
\[
\varphi_\star \colon PrB_\nn \longrightarrow \PLB\nn.
\]
Indeed, two pure geometric braids $\beta^\prime$ and $\beta^{\prime\prime}$ are equivalent if and only if there is an ambient isotopy of $\Rr^3 \times [0, 1]$ from the identity map to a self-homeomorphism that maps $\beta^\prime$ to $\beta^{\prime\prime}$. That by construction would be an isotopy (so in particular a homotopy) between the two associated loops in~$\PLB\nn$. Moreover products are preserved, so $\varphi_\star$ is a isomorphism.
\endproof

\begin{thm}
\label{T:flying}
For $\nn \geq 1$, there is an isomorphism between the ribbon braid group $\rB\nn$ and the loop braid group~$\LB\nn$.
\end{thm}
\proof
We recall that the untwisted ring group $\UR\nn$ (Definition~\ref{D:configurations}) is the fundamental group of $\UnRi\nn$, which is the quotient of $\PUnRi\nn$ by the symmetric group on $\nn$ components. As in Proposition~\ref{P:pure} we fix an element $\beta=\big( C_1(t), \dots, C_\nn(t)\big)$ of $\rB\nn$ and define a map:
\[
\hat{\varphi}(\beta) \colon [0, 1] \longrightarrow \LBE\nn
\]
by $t \mapsto \big[ C_1(t), \dots, C_\nn(t)\big]$. The element $\hat{\varphi}(\beta)$ is a loop in the configuration space $\UR\nn$. This loop corresponds to an element of $\LB\nn$ through the isomorphism from Proposition~\ref{P:PureUnrestricted}. Then $\hat{\varphi}$ induces an homomorphism 
\[
\hat{\varphi}_\star \colon \rB\nn \longrightarrow \LB\nn.
\]

We consider the following diagram:
\[
\begin{CD}
1 @>>> PrB_\nn @>>>  rB_\nn @>>> S_\nn @>>>1 \\
   @.  @V{\cong}V{\varphi_\star}V          @VV{\hat{\varphi}_\star}V        @|  \\
1 @>>> \PLB\nn @>>> \LB\nn @>>> S_\nn @>>> 1 .
\end{CD}
\]
It is commutative by construction of $\varphi$ and~$\hat{\varphi}$. By applying the five lemma, the statement is proved.
\endproof

\subsection{Extended ribbon braids}
We introduce a new topological object, that is the topological realization of extended loop braids. 

\begin{defn}
Let $\nn \geq 1$. An \emph{extended geometric ribbon braid} is a geometric ribbon braid, with condition \eqref{wen} from Definition \ref{D:ribbonBraid} removed, and condition~\eqref{bordo} replaced by:
\begin{enumerate}
\item[\eqref{bordo}$^{\prime \prime}$] the boundary of each annulus $\partial A_i$ is a disjoint union~$C_i \sqcup C_j$, for $C_i \in \partial_0 B^4$ and for some~$C_j \in \partial_1 B^4$. 
\end{enumerate}

The group of \emph{extended ribbon braids}, denoted by $\rBE\nn$, is the group of equivalence classes of geometric ribbon braids up to continuous deformations through the class of extended geometric ribbon braids fixing the boundary circles, equipped with the natural product given by stacking and reparametrizing. The unit element for this product is the \emph{trivial ribbon braid}~$U= \bigsqcup_{i\in \{1, \dots, n \}} C_i \times [0, 1]$.
\end{defn}

\begin{defn}
The group of \emph{extended pure ribbon braids}, denoted by $PrB_n^{ext}$, is the group of equivalence classes of pure geometric extended ribbon braids, defined as pure geometric ribbon braids (Definition \ref{D:pureRibbon}) with condition \eqref{wen} removed and condition \eqref{bordo}$^\prime$ replaced by
\begin{enumerate}
\item[\eqref{bordo}$^{\prime \prime \prime}$] $\partial A_i = C_i \times \{0, 1\} $ for all~$i \in \{1, \dots, n\}$ .
\end{enumerate} 
\end{defn}

\begin{rmk}
\label{R:extending}
Theorem~\ref{T:isotopy} extends to extended ribbon braid groups $\rBE\nn$: in fact the proof uses the evaluation filtration $\tilde{\ep}\colon \Homeo(B^3) \to \mathcal{L}_{\nn}$. We recall that, for $\nn \geq 1$ fixed,  $\mathcal{L}_{\nn}$ is the space of configurations of $\nn$ unordered, disjoint, unlinked trivial links in~$B^3$, that are not restrained of parallel planes, \ie, can make $180^{\circ}$ flips. 
\end{rmk}
\begin{rmk}
The monotony condition allows us to consider the interval $I$ in $B^4 = B^3 \times I$ as a time parameter, and to think of an extended ribbon braid as a trajectory $\beta=\big( C_1(t), \dots, C_\nn(t)\big)$ of circles in $B^3 \times I$. This trajectory corresponds to a parametrization of the extended ribbon braid. That corresponds to an interpretation in terms of flying rings allowed to ``flip''. 
\end{rmk}

The following statements, proved for the case of ribbon braids, hold also in the extended case.
\begin{prop} 
\label{P:Extpure}
For $\nn \geq 1$, there is an isomorphism between the pure extended ribbon braid group $PrB_\nn^{ext}$ and the pure extended loop braid group~$\PLBE\nn$.
\end{prop}
\proof 
To prove this statement it is enough to follow the proof of  Proposition~\ref{P:pure}, replacing $\PUR\nn$ with~$\PR\nn$, and considering the isomorphism between and $\PR\nn$ and~$\PLBE\nn$.
\endproof

\begin{thm}
\label{T:Extflying}
For $\nn \geq 1$, there is an isomorphism between the extended ribbon braid group $\rBE\nn$ and the extended loop braid group~$\LBE\nn$ .
\end{thm}
\proof
To prove this statement it is enough to follow the proof of  Theorem~\ref{T:flying}, replacing $\UR\nn$ with $\R\nn$, and considering the isomorphism between $\R\nn$ and~$\LBE\nn$.
\endproof


\section{Welded diagrams and broken surfaces}
\label{S:welded}
There are two kind of projections of loop braids, seen as ribbon braids: the first one is a $2$-dimensional diagrammatical representation, while the second one will be a representation through $3$-dimensional surfaces. 
 
\subsection{Welded diagrams}
\label{SS:WeldedBraidDiagrams}
\begin{defn}
\label{D:strand}
A \emph{strand diagram on $\nn$ strings} is a set of oriented arcs in $\Rr^2$, monotone with respect to the second coordinate, from the points $(0, 1), \dots, (0, \nn)$ to $(1, 1), \dots, (1, \nn)$. The arcs are allowed to have double points of three kinds, called \emph{classical positive}, \emph{classical negative} and \emph{welded} as in Figure~\ref{F:Crossings}.

\begin{figure}[hbt]
\centering
\includegraphics[scale=0.6]{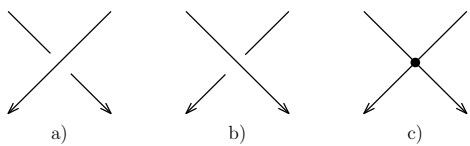}
\caption{a) Classical positive crossing, b) Classical negative crossing, c) Welded crossing.}
\label{F:Crossings}
\end{figure}

Let us assume that the double points occur at different $x$-coordinates. Then a strand diagram determines a word in the elementary diagrams illustrated in Figure~\ref{F:Crossings}. We call $\sig\ii$ the elementary diagram representing the $(i+1)$-th strand passing over the $i$-th strand, and $\rr\ii$ the welded crossing of the strands $i$ and~$(i+1)$. The set of strand diagrams on $\nn$ strings will be denoted by~$\mathcal{D}_\nn$.
\begin{figure}[hbt]
\centering
\includegraphics[scale=0.6]{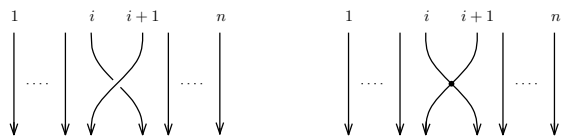}
\caption{The elementary diagrams $\sig\ii$ and $\rr\ii$.}
\label{F:SigRho}
\end{figure}
\end{defn}

\begin{figure}[hbt]
\centering
\includegraphics[scale=.6]{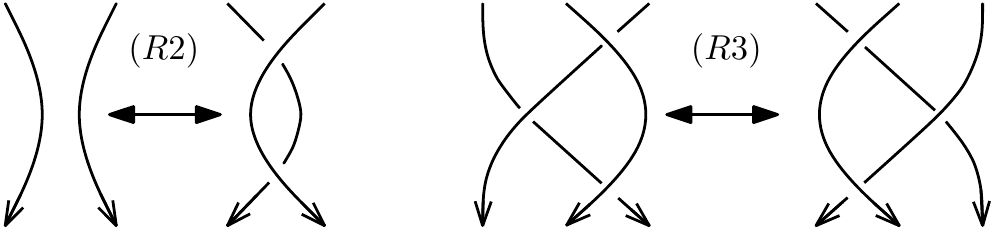}
\caption{Classical Reidemeister moves for braid-like objects.}
\label{F:Classical}
\end{figure}

\begin{figure}[hbt]
\centering
\includegraphics[scale=.6]{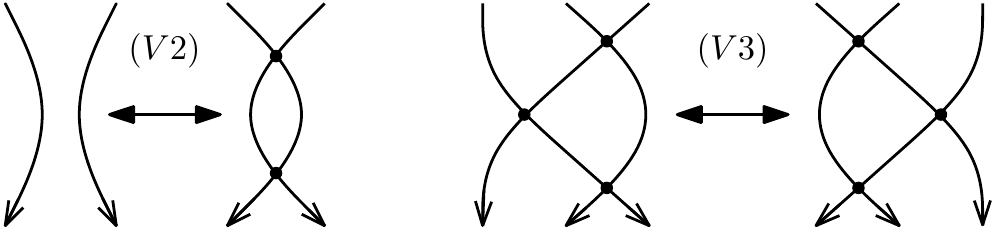}
\caption{Virtual Reidemeister moves for braid-like objects.}
\label{F:Virtual}
\end{figure}

\begin{figure}[hbt]
\centering
\includegraphics[scale=.6]{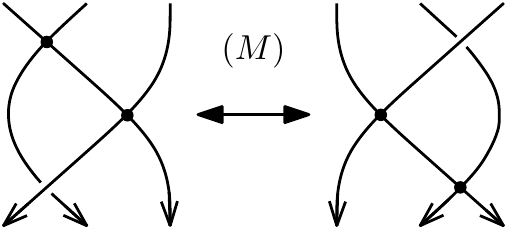}
\caption{Mixed Reidemeister move.}
\label{F:Mixed}
\end{figure}

\begin{figure}[hbt]
\centering
\includegraphics[scale=.6]{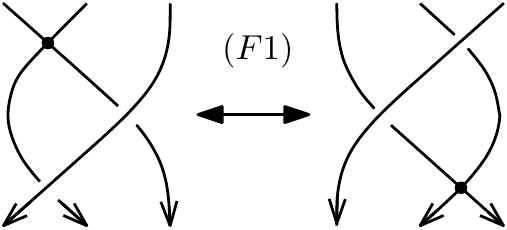}
\caption{Welded Reidemeister move.}
\label{F:Welded}
\end{figure}

\begin{defn}
\label{D:bWelded}
A \emph{welded braid} is an equivalence class of strand diagrams under the equivalence relation given by planar isotopy and the following moves: 
\begin{itemize}
\item Reidemester moves: Figure~\ref{F:Classical};
\item virtual Reidemeister moves: Figure~\ref{F:Virtual};
\item mixed Reidemeister moves: Figure~\ref{F:Mixed};
\item welded Reidemeister moves: Figure~\ref{F:Welded}.
\end{itemize} 
This equivalence relation is called \emph{welded Reidemeister equivalence}. We denote classes by representatives.
For $\nn \geq 1$, the \emph{group of welded braids} or \emph{welded braid group} on $\nn$ strands, denoted by~$\WB\nn$ is the group of equivalence classes of strand diagrams by welded Reidemeister equivalence. The group structure on these objects is given by: stacking and rescaling as product, mirror image as inverse, and the trivial diagram as identity.
\end{defn}

\begin{rmk}
\label{R:conventions}
We adopt the convention of reading the diagrams from top to bottom, and the corresponding braid words from left to right. Welded braids have been introduced in~\cite{FRR} as a graphic representation of presentation~\eqref{E:Upresentation}, given in their paper as a presentation of the braid-permutation groups~$\BP\nn$. 
Then notation used here for elementary diagrams does not carry any ambiguity with the notation used for the presentation~\eqref{E:Upresentation} of $\UR\nn$ in Section~\ref{S:configuration}, and for the automorphisms of the free group in Section~\ref{S:automorphisms}. 
\end{rmk}

\subsection{Broken surfaces}
\label{SS:BrokenBraidSurfaces}
\begin{figure}[hbt]
\centering
\includegraphics[scale=.5]{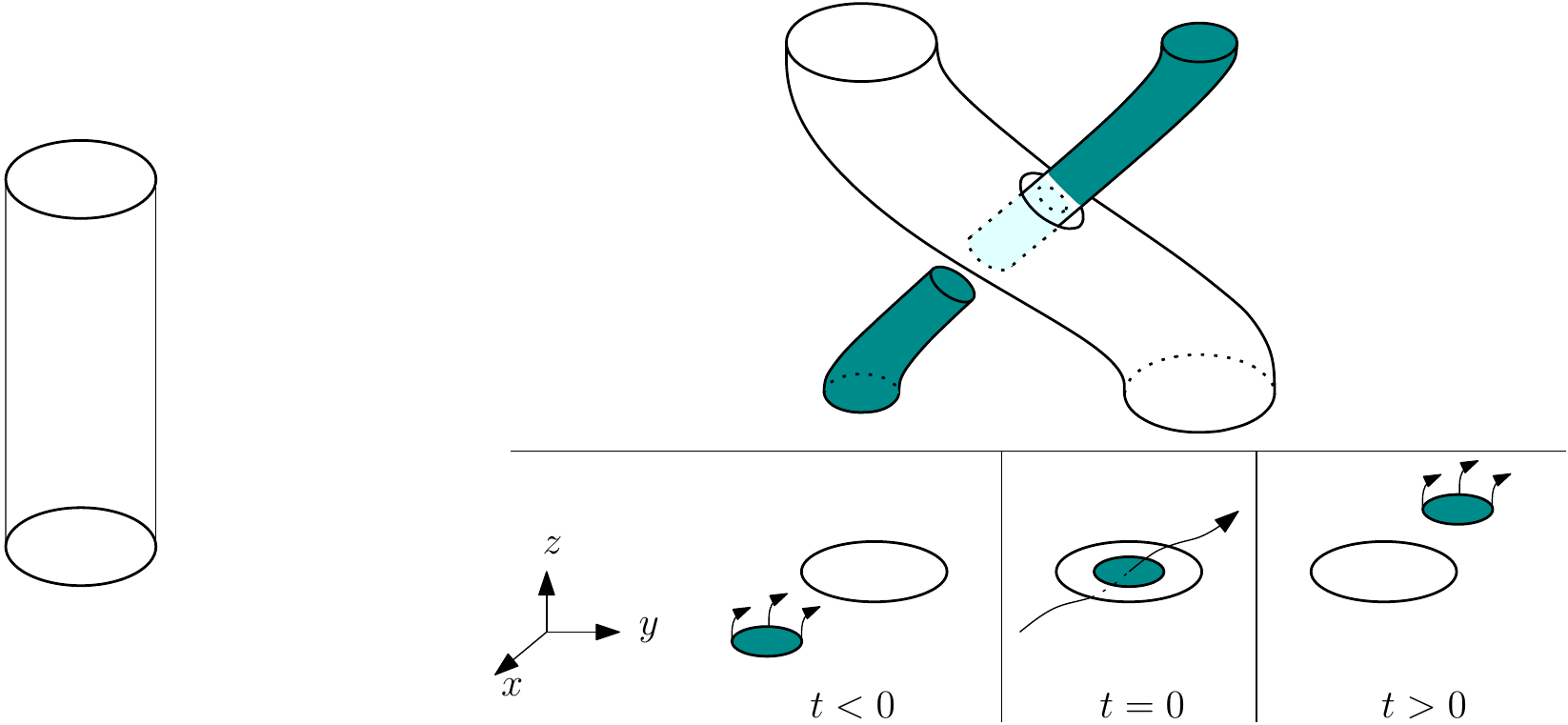}
\caption{Symmetric broken surfaces are locally homeomorphic to a cylinder, or to a crossing with decoration that indicates the order on the preimages of the singularities.}
\label{F:BrokenDeco}
\end{figure}

A diagram of a classical braid is a projection in general position of the braid on the plane, with crossing information specified by deleting a neighbourhood of the underpassing arcs. We introduce similar diagrams for geometric ribbon braids, which are surfaces in the $4$-dimensional space: these diagrams are projections in general position of ribbon braids in the $3$-dimensional space, and are called \emph{broken surface diagrams}. 
This representation has first been introduced for the group of motions of a collection of $n$ unknotted, unlinked, oriented circle in~\cite{BWC}. They have defined broken surfaces adapting a drawing style from Carter and Saito's work on surfaces in dimension $4$~\cite{CS}. We use notations introduced at the beginning of Section~\ref{S:ribbon}.
\begin{defn}
\label{D:broken}

Let $A_1, \dots, A_\nn$  be locally flat embeddings in~$\stackrel{*}{B^3}$ of $\nn$ disjoint copies of the oriented annulus~$S^1 \times I$.  We say that  
\[S = \bigcup_{i\in \{1, \dots, n \}} A_i 
\]
is a \emph{braid broken surface diagram} if:
\begin{enumerate}
\item for each $i \in \{1, \dots, \nn\}$, the oriented boundary $\partial A_i$ is the disjoint union $C_i \sqcup C_j$, for $C_i$ in $\partial_0 B^3$ and for some $C_j$ in $\partial_1 B^3$. The orientation induced by $A_i$ on $\partial A_i$ coincides with the orientation of one of the two boundary circles $C_i$ and~$C_j$;
\item  it is transverse to the lamination $\bigcup_{t \in I} B^2 \times \{t\}$ of~$B^3$, that is: at each parameter~$t$, the intersection between $S$ and $B^2 \times{t}$ is a collection of exactly $\nn$ circles, not necessarily disjoint;
\item the set of connected components of singular points in $S$, denoted by $\Sigma(S)$, consists of flatly transverse disjoint circles in $(\cup_{i = 1}^n \mathring{A_i})$.
\end{enumerate}
Moreover, for each element of~$\Sigma(S)$, a local ordering is given on the two circle preimages. By convention this ordering is specified on the diagram by erasing a small neighbourhood of the lower preimage in the interior of the annulus it belongs to. Note that this is the same convention used for classical braid diagrams. 
Moreover a broken surface diagram is said to be \emph{symmetric} if it is locally homeomorphic to the surfaces in Figure~\ref{F:BrokenDeco}, which means that, for each pair of preimage circles, the following properties are satisfied:
\begin{enumerate}
\item one of the preimage circles is essential in $\bigcup_{i = 1}^n \mathring{A_i}$ and the other is not;
\item there is a pairing of the elements of $\Sigma(S)= \bigsqcup_r \{c^r_1, c^r_2\}$ such that, for each~$r$, the essential preimages of~$c^r_1$ and~$c^r_2$
\begin{enumerate}
\item are respectively lower and higher than their non essential counterparts with respect to the associated order ;
\item bound an annulus in $\bigcup_{i=1}^n \mathring{A_i}$;
\item this annulus avoids $\Sigma(S)$.
\end{enumerate}
\end{enumerate}
\end{defn}

Let $b$ be a ribbon braid, and consider a projection $B^4 \to B^3$ in general position of~$b$: the following result allows us to consider braid broken surface diagrams as $3$-dimensional representations of ribbon braids.

\begin{lem}
\label{L:RibbonEBroken}
Any generic projection of a ribbon braid from $B^4$ to $B^3$ is a braid broken surface diagram. Conversely any braid broken surface diagram is the projection of a ribbon braid. 
\end{lem}
\proof
This statement is proved in \cite{YAN} for locally flat embedding of $2$-spheres in $\Rr^4$. In \cite{ABMW} it is noted that the arguments, which are local, apply to the case of ribbon tubes. Since ribbon tubes and ribbon braids locally behave the same way, the arguments of \cite{YAN} applies to our case. 
\endproof

Moreover, we have the following result:

\begin{lem}
\label{L:symmetricBroken}
Any ribbon braid can be represented by a symmetric broken surface diagram.
\end{lem}
\proof
This statement is proved in~\cite[Lemma~2.13]{ABMW} for ribbon tubes, adapting results from \cite{KS} and \cite{YAJ}. The proof for ribbon tubes adapts without modifications to ribbon braids.
\endproof

\begin{figure}[htb]
	\centering
		\includegraphics[scale=0.5]{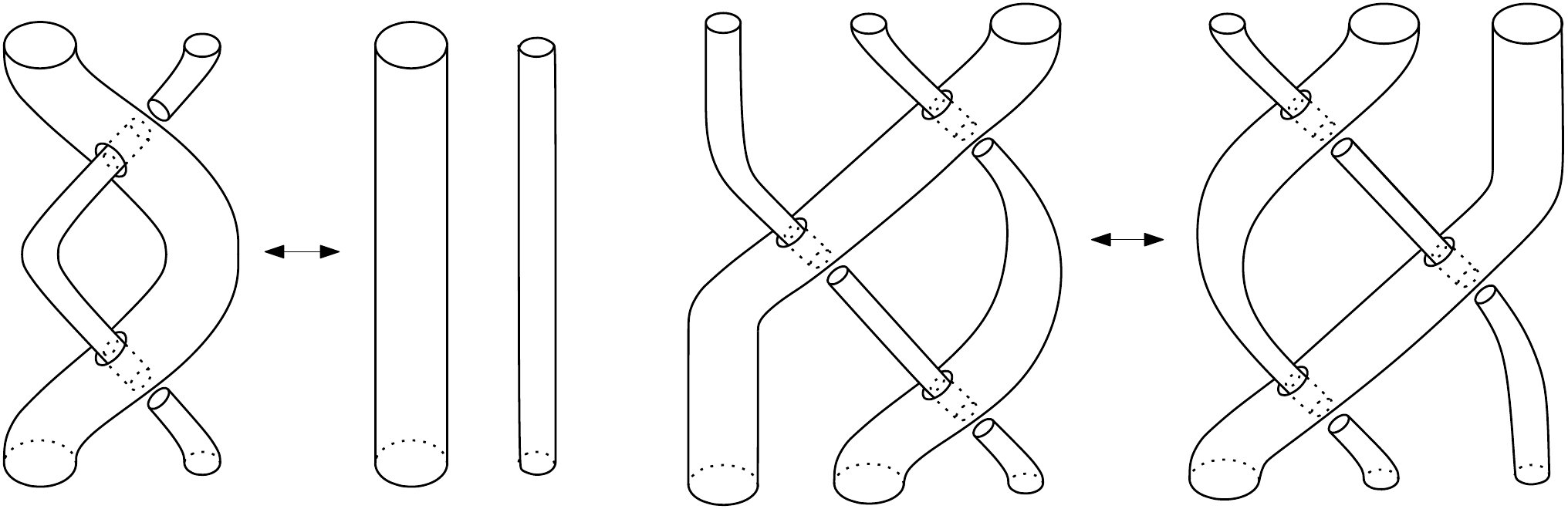} 
	\caption{Broken Reidemeister moves.}
	\label{F:BrokenReid}
\end{figure}

\begin{rmk}
If two symmetric braid broken surface diagrams differ by one of the ``broken Reidemeister moves'' in Figure~\ref{F:BrokenReid}, then the associated ribbon braids are isotopic~\cite[Remark~2.15]{ABMW}.
\end{rmk}

\begin{figure}[htb]
	\centering
		\includegraphics[scale=0.7]{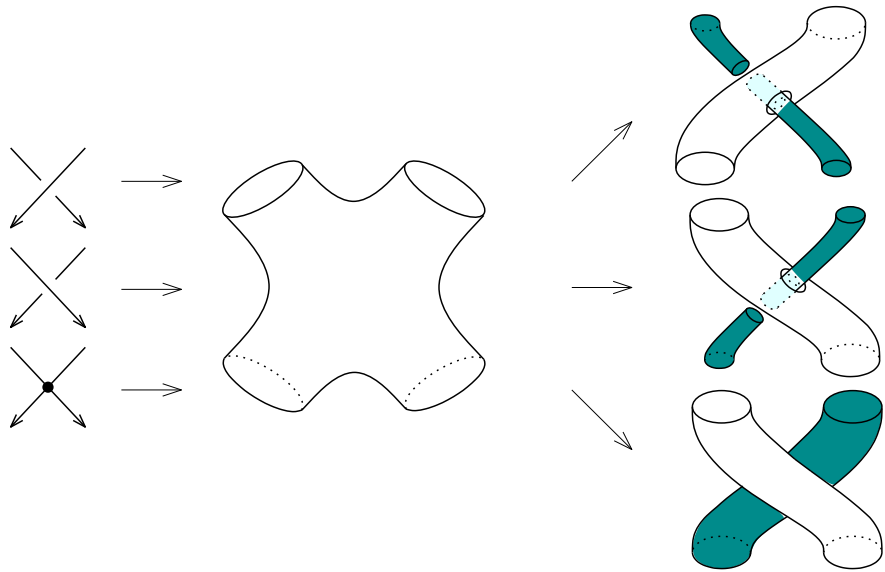} 
	\caption{Punctured sphere associated to a crossing of the welded braid.}
	\label{F:tubemap}
\end{figure}

Passing through symmetric braid broken surfaces, $4$-dimensional ribbon braids can be described using $2$-dimensional welded braids. Let $b$ be a welded braid. We associate to it a symmetric braid broken surface diagram in the following way (see for details \cite{SAT} and~\cite{YAJ2}). 
Consider $B^2$ and embed it as $B^2 \times \{\frac{1}{2}\}$ into $B^3$. Taken a tubular neighbourhood $N(b)$ of~$b$, we have that $\partial_{\ep}N(b) = \sqcup_{i \in \{1, \dots, n\}} D_i \times \ep_i$ where~$\ep_i \in \{0, 1\}$. Each crossing is sent to a $4$-punctured sphere. Then, according to the partial order defined on the double points of welded braid diagrams, we modify the sphere into the broken surfaces shown in Figure~\ref{F:tubemap}).

\begin{defn}
\label{D:tube}
We define a map $Tube \colon \WB\nn \to \rB\nn$ that associates to a welded braid~$b$, the ribbon braid associated to the symmetric broken surface resulting from the preceding construction.
\end{defn}

For general ribbon knotted objects, such as ribbon tangles~\cite{DF}, ribbon tubes~\cite{ABMW} and knotted spheres~\cite{YAN}, the Tube map is well defined and surjective.  For these more general objects, it is still unproved if the Tube map is injective~\cite{Bn}. However, Brendle and Hatcher proved that $\UR\nn$ is isomorphic to the presented group \eqref{E:Upresentation}, represented by $\WB\nn$ \cite{BH}. Through the isomorphism between $\UR\nn$ and $\LB\nn$ (Proposition~\ref{P:PureUnrestricted}) and the isomorphism between $\LB\nn$ and $\rB\nn$ (Theorem~\ref{T:flying}), we have the following result, that gives us a graphical interpretation of loop braids.
\begin{thm}
\label{T:IsoTube}
The map $Tube \colon \WB\nn \to \rB\nn$ is an isomorphism.
\end{thm}

Let $\widetilde{Tube}$ the composition of $Tube$ with the isomorphism between $\rB\nn$ and $\LB\nn$. We have:

\begin{cor}
\label{C:Tube}
The map $\widetilde{Tube} \colon \WB\nn \to \LB\nn$ is an isomorphism.
\end{cor}

\subsection{Extended welded diagrams and broken surfaces}

\begin{figure}[!htb]
\centering
\includegraphics[scale=0.6]{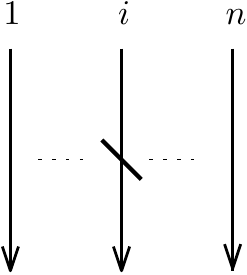}
\caption{The elementary diagram~$\tau_\ii$.}
\label{F:Tau}
\end{figure}

\begin{figure}[htb]
\centering
\includegraphics[scale=0.6]{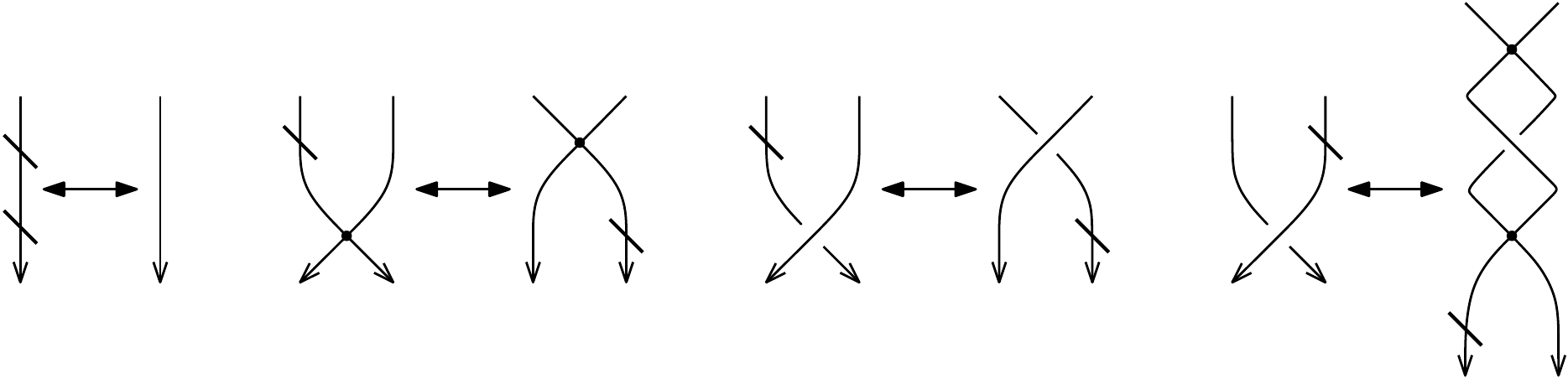}
\caption{Extended Reidemeister moves for braid-like objects.}
\label{F:ReidWen}
\end{figure}

We introduce new kind of local elementary diagrams, called~$\tau_\ii$, for $I=1, \dots, \nn$ as in Figure~\ref{F:Tau}, 

\begin{defn}
\label{D:bExteWelded}
An \emph{extended welded braid} is an equivalence class of strand diagrams under the equivalence relation given by planar isotopy and the following moves: 
\begin{itemize}
\item usual Reidemeister moves: Figure~\ref{F:Classical};
\item virtual Reidemeister moves: Figure~\ref{F:Virtual};
\item mixed Reidemeister moves: Figure~\ref{F:Mixed};
\item welded Reidemeister moves: Figure~\ref{F:Welded};
\item extended Reidemester moves: Figure~\ref{F:ReidWen}.
\end{itemize} 
This equivalence relation is called \emph{extended welded Reidemeister equivalence}. We denote classes by representatives.
The \emph{group of extended welded braids} or \emph{extended welded braid group}, denoted by~$\WBE\nn$ is the group of equivalence classes of strand diagrams by extended welded Reidemeister equivalence. The group structure on these objects is given by: stacking and rescaling as product, mirror image as inverse, and the trivial diagram as identity.
\end{defn}

\begin{figure}[htb]
\centering
\includegraphics[scale=0.6]{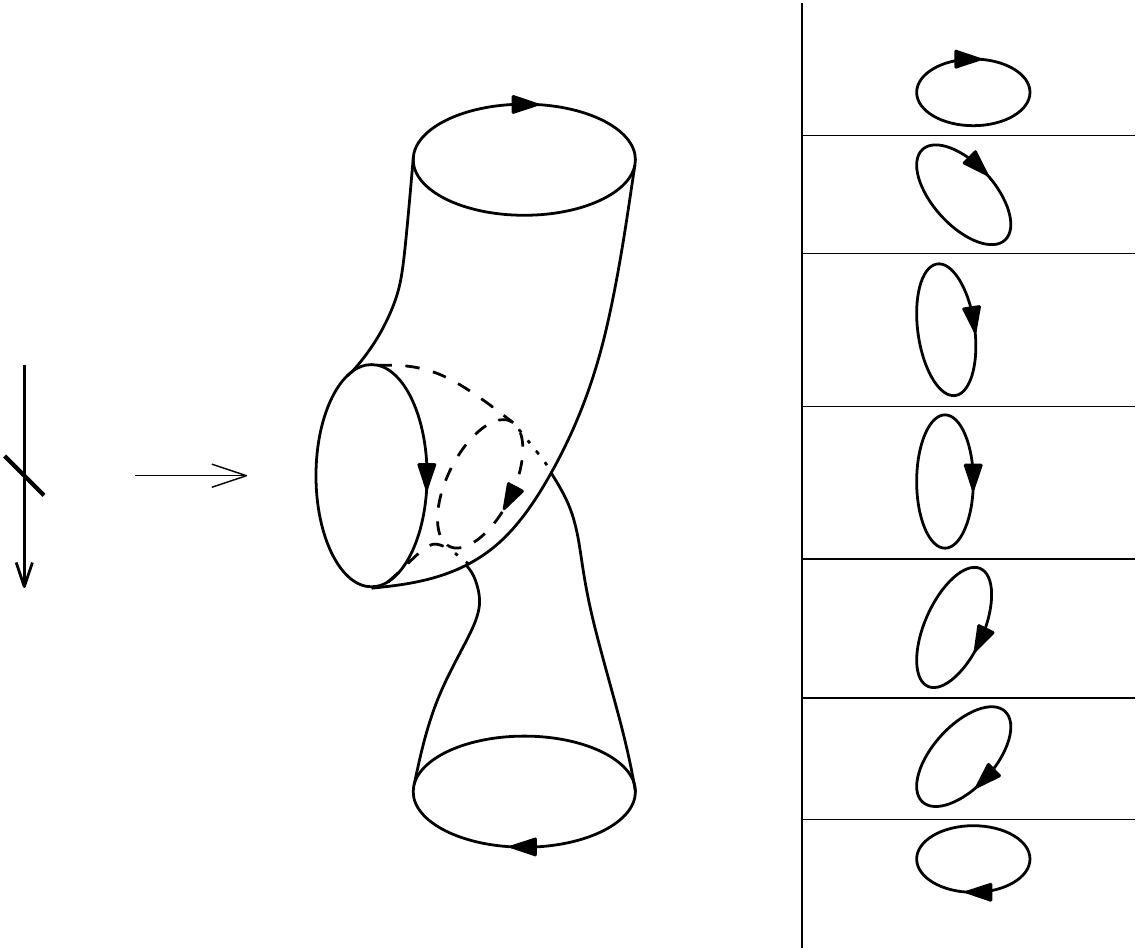}
\caption{The behaviour of the Tube map when applied to wens.}
\label{F:TubeWen}
\end{figure}

The $Tube$ map can be extended to $\WBE\nn$, associating to elementary diagrams $\tau_i$ a tube with a wen (Figure~\ref{F:TubeWen}). Brendle and Hatcher proved that $\R\nn$ is isomorphic to the presented group \eqref{E:Rpresentation}, this second one represented by~$\WBE\nn$~\cite{BH}. Through the isomorphism between $\R\nn$ and $\LBE\nn$ (Theorem~\ref{T:PureRing}),  and the isomorphism between $\LBE\nn$ and $\rBE\nn$ (Theorem~\ref{T:Extflying}), we have the following result, that gives us a graphical interpretation of extended loop braids.

\begin{thm}
The map $Tube \colon \WBE\nn \to \rBE\nn$ is an isomorphism.
\end{thm}

\section{Gauss diagrams}
\label{S:Gauss}
In this section we give a combinatorial description of loop braids through Gauss diagrams. 

\begin{figure}[htb]
\centering
\includegraphics[scale=0.7]{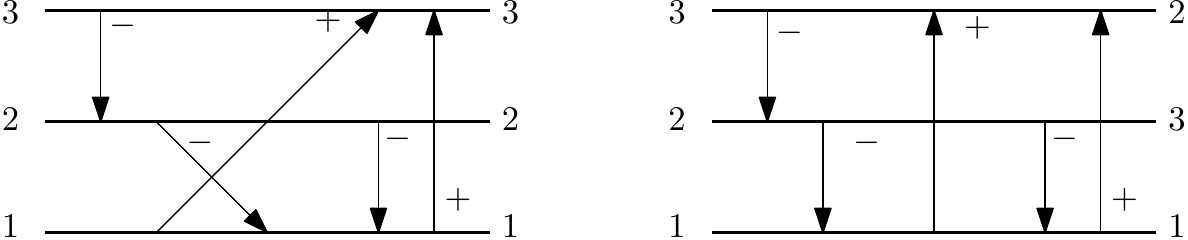}
\caption{Examples of Gauss diagrams.}
\label{F:Gauss1}
\end{figure}

\begin{defn}Let $\nn \geq 1 $. A \emph{Gauss diagram} is a set of signed and oriented arrows on $n$ ordered and oriented horizontal intervals, together with a permutation $\sigma \in S_n$ (see Figure~\ref{F:Gauss1}). The endpoints of the arrows are divided in two sets: the set of \emph{heads} and the set of \emph{tails}, defined by the orientation of the arrow. The right extremity of the $i$th horizontal interval is labelled with~$\sigma(i)$. 
\end{defn}

To make the link between Gauss diagrams and the group of welded braids, we introduce the group of virtual braids. This allows us to see the welded braid group as a quotient of virtual braids. Generally speaking, Gauss diagrams turn out to be a useful tool to investigate properties of various remarkable quotients of the virtual braid group (see for example \cite{ABMW2, ABMW3}). Recall Definition~\ref{D:strand} of strand diagrams.

\begin{defn}
Two strand diagrams on $\nn$ strands are \emph{virtual equivalent} if they are related by planar isotopy and a finite number of the following moves:
\begin{itemize}
\item virtual Reidemeister moves: $(V2)$ and $(V3)$,  Figure~\ref{F:Virtual};
\item mixed Reidemeister moves: $(M)$, Figure~\ref{F:Mixed}.
\end{itemize}
Two strand diagrams on $\nn$ strands are \emph{welded equivalent} if they are related by planar isotopy and a finite number of the following moves:
\begin{itemize}
\item virtual Reidemeister moves: $(V2)$ and $(V3)$,  Figure~\ref{F:Virtual};
\item mixed Reidemeister moves: $(M)$, Figure~\ref{F:Mixed};
\item welded Reidemeister moves: $(F1)$, Figure~\ref{F:Welded}.
\end{itemize}
Two strand diagrams on $\nn$ strands  are \emph{virtual Reidemeister equivalent} if they are related by planar isotopy and a finite number of the following moves:
\begin{itemize}
\item Reidemester moves: $(R2)$ and $(R3)$, Figure~\ref{F:Classical};
\item virtual Reidemeister moves: $(V2)$ and $(V3)$,  Figure~\ref{F:Virtual};
\item mixed Reidemeister moves: $(M)$, Figure~\ref{F:Mixed}.
\end{itemize}
\end{defn}
We recall that adding welded Reidemeister moves to virtual Reidemeister equivalence, one obtains welded Reidemeister equivalence (Definition~\ref{D:bWelded}).

\begin{defn}
\label{D:VirtualBraidGroup}
For $\nn \geq 1$, the \emph{virtual braid group} $\VB\nn$ is the group of equivalence classes of strand diagrams on $n$ strands with respect to the virtual Reidemeister equivalence. We call \emph{virtual braid diagram} an element of this group.
\end{defn}

Other equivalent definitions of virtual braid groups have been introduced for instance in~\cite{Ver01, Bn,Kam}. 


\begin{rmk}
\label{R:WeldedQuotient}
For $\nn \geq 1$, the welded braid group $\WB\nn$ is a quotient of the virtual braid group $\VB\nn$ under the relation given by moves of type~$(F1)$. 
\end{rmk}

To every virtual braid diagram $\beta$ we can associate a Gauss diagram with the following construction. Let $\beta$ be a strand diagram on $n$ strands. The \emph{Gauss diagram} associated to~$\beta$, denoted by~$G(\beta)$, is a Gauss diagram on $\nn$ intervals satisfying the following properties:
\begin{enumerate}
\item for each strand of $\beta$ there is an associated interval of~$G(\beta)$;
\item the endpoints of the arrows of $G(\beta)$ correspond to the preimages of the classical (positive and negative) crossings of~$\beta$;
\item the order of the endpoints of the arrows on an interval of $G(\beta)$ correspond to the order that the preimages associated to the endpoints on the strand of $\beta$;
\item the arrows are pointing from the overpassing strand to the underpassing strand, when considering the usual convention on strand diagrams of drawing a break on the underpassing strand;
\item the permutation of $G(\beta)$ corresponds to the permutation defined by~$\beta$.
\end{enumerate}

\begin{rmk}
\label{R:reparametrization}
We remark that the Gauss diagram associated to a strand diagram has pairwise distinct arrows, and each arrow connects different intervals. By reparametrizing the intervals, the arrows can always be drawn to be vertical and at different $t$ coordinates, for $t \in I$, where $t$ is the horizontal coordinate.
\end{rmk}

\begin{defn}
Gauss diagrams respecting the conditions given in Remark~\ref{R:reparametrization} are called \emph{braid Gauss diagrams}. Their set will be denoted by~$bGD_n$.
\end{defn}

\begin{defn}
Gauss diagrams respecting the conditions given in Remark~\ref{R:reparametrization} are called \emph{braid Gauss diagrams}. Their set will be denoted by~$bGD_n$.
\end{defn}

\begin{figure}[!htb]
\centering
\includegraphics[scale=.6]{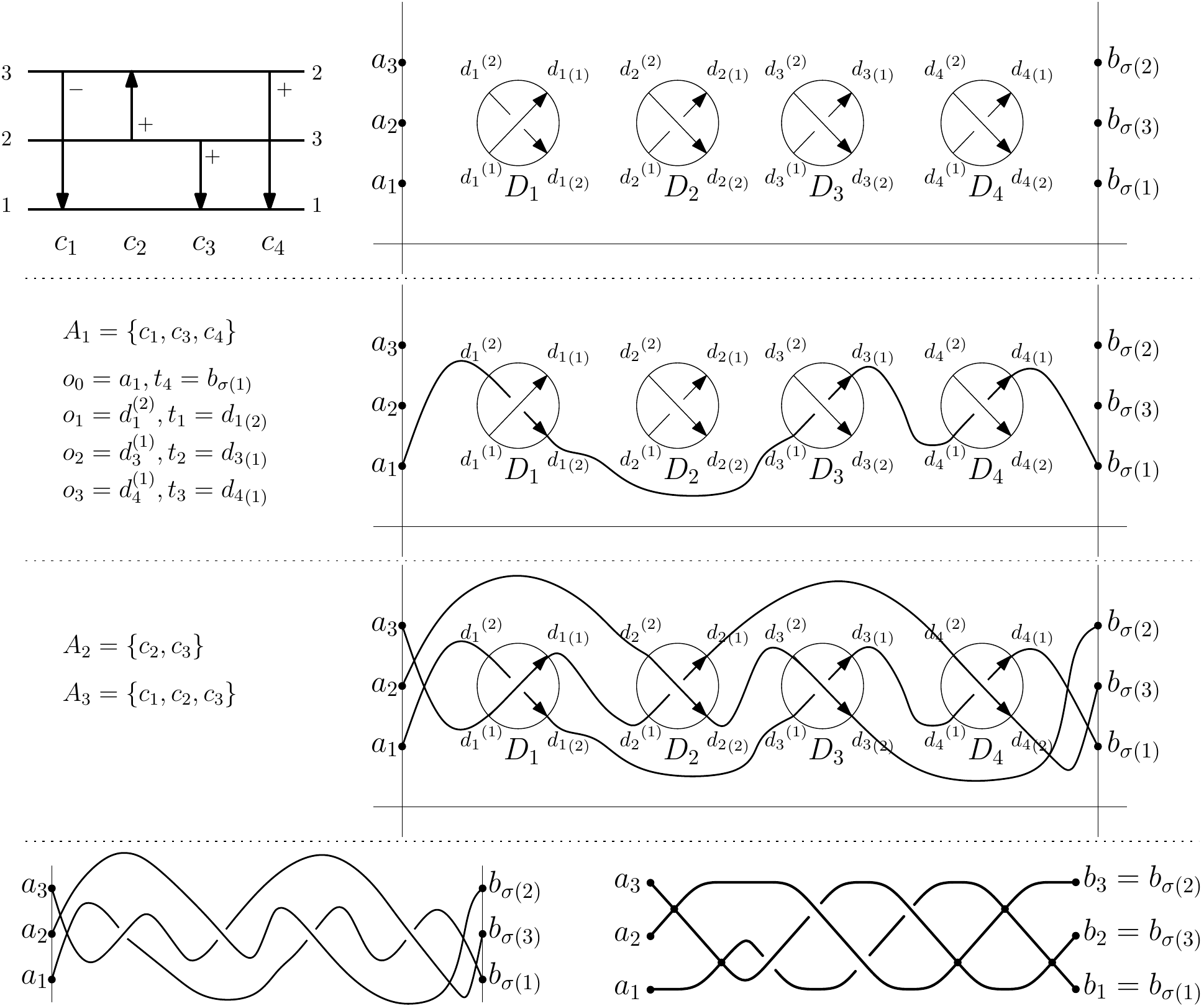}
\caption{From a braid Gauss diagram to a strand diagram.}
\label{F:GaussEsempio}
\end{figure}

The following result is proved by Cisneros in~\cite{Ci}. We give an example in Figure~\ref{F:GaussEsempio}.
\begin{thm}[{ \cite[Theorem~2.10]{Ci}}]
\label{T:Bruno}
The following statements hold:
\begin{enumerate}[label=(\roman*)]
\item \label{i:1} For each braid Gauss diagram $g$ there exists a braid $\beta \in \mathcal{D}_\nn$ such that $G(\beta)=g$. 
\item \label{i:2} Given two virtual braids $\beta_1$ and $\beta_2$ in $\VB\nn$, $G(\beta_1)=G(\beta_2)$ if and only if $\beta_1$ and $\beta_2$ are \emph{virtual equivalent}.
\end{enumerate}
\end{thm}
%
%
%
%

Defining equivalence relations on~$bGD_n$, Cisneros proves the existence of a bijection between a quotient of $bGD_n$ and $\VB\nn$~\cite[Proposition~2.24]{Ci}. This can be extended to~$\WB\nn$, giving us the last isomorphism of this paper. 

\begin{defn} Let $g$ and $g^\prime$ be two Gauss diagrams. A \emph{Gauss embedding} is an embedding $\phi \colon g^\prime \to g$ that sends each interval of $g^\prime$ to a subinterval of~$g$, and which sends each arrow of $g^\prime$ to an arrow of $g$ respecting  orientation and  sign. Note that there is no condition on the permutations associated to  $g$ and~$g^\prime$.
\end{defn}

\begin{figure}[!htb]
\centering
\includegraphics[scale=.7]{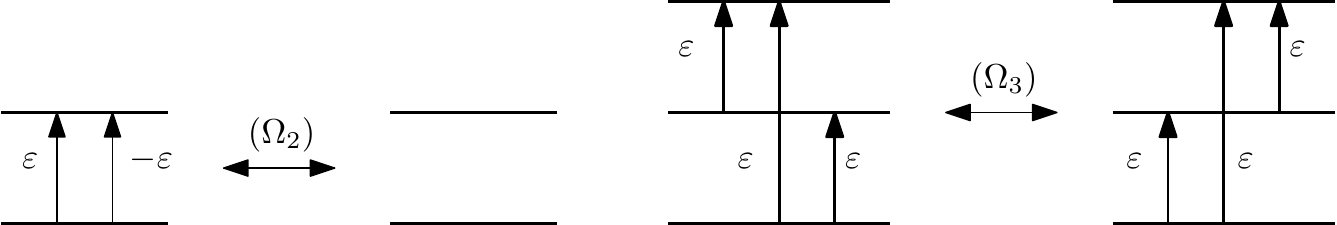}
\caption{Equivalence moves on braid Gauss diagrams, where the $\ep$ is the sign of the arrows.}
\label{F:GaussReid}
\end{figure}

\begin{figure}[!htb]
\centering
\includegraphics[scale=.7]{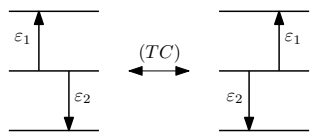}
\caption{Tail Commute move where the $\ep_i$ are the signs of the arrows.}
\label{F:GaussTail}
\end{figure}

\begin{defn}
\label{D:GaussReid}
Two braid Gauss diagram are \emph{Reidemeister equivalent} if they are related by a finite sequence of moves $(\Omega_2)$ and $(\Omega_3)$ as in Figure~\ref{F:GaussReid}. They are \emph{wReidemeister equivalent} if they are related by~$(\Omega_2)$,~$(\Omega_3)$, and Tail Commute moves $(TC)$ as in Figure~\ref{F:GaussTail}. For $\nn \geq 1$, the group of \emph{welded Gauss diagrams}, denoted by~$wG_n$, the group of equivalence classes of $bGD_n$ with respect to wReidemeister equivalence. 
\end{defn}

\begin{rmk}
Performing one of the moves described in Definition~\ref{D:GaussReid} on a braid Gauss diagram $g$ means choosing a Gauss embedding in $g$ of the braid Gauss diagram depicted on one hand of the equivalences in Figures~\ref{F:GaussReid} and~\ref{F:GaussTail}, and replacing it with the braid Gauss diagram on the other hand. Since there is no condition on the permutation associated, when embedding a Gauss diagram, this means that in performing the moves the strands can be vertically permuted.
\end{rmk}

The following result gives us an interpretation of $\LB\nn$ in terms of Gauss diagrams, recalling that $\LB\nn$ is isomorphic to~$\WB\nn$.

\begin{thm}
\label{T:Gauss}
For $\nn \geq 1$, there is a bijective correspondence between the group of welded Gauss diagrams $wG_n$ and the welded braid group~$\WB\nn$.
\end{thm}
\proof
By Theorem~\ref{T:Bruno} we know  that there is a bijective correspondence between virtual equivalent braids and braid Gauss diagrams. Therefore we need to prove that if two strand diagrams are related by a classical Reidemeister move~$(R2)$,~$(R3)$, or a welded Reidemeister move~$(F1)$, then their braid Gauss diagram are wReidemeister equivalent via moves~$(\Omega_2)$,~$(\Omega_3)$ and~$(TC)$, and viceversa.  

Moves $(\Omega_1)$ and $(\Omega_2)$ are treated in~\cite[Proposition~2.24]{Ci}. Let $\beta$ and $\beta^\prime$ be two welded braids that differ by an $(F1)$ move, and suppose that the strands involved are~$a$, $b$ and $c$ in~$\{1, \dots, n\}$. 
Then, up to isotopy, we can deform the two braids so that they coincide outside of a portion of the braid diagram which only contains the involved crossings. In this portion the diagrams look like in Figure~\ref{F:GaussReid}, and the Gauss diagrams associated differ by a $(TC)$ move. 

Let $g$ and $g^\prime$ be two Gauss diagrams and let $a, b$ and $c$ in $\{1, \dots, n\}$ be pairwise different. Assume they differ by a $(TC)$ move. Then there exists a sub-portion of the diagram that contains only the arrows involved in the $(TC)$ move. For Theorem~\ref{T:Bruno} there exists a strand diagram $\beta$ that looks like the right or the left handside of move $(F1)$, in Figure~\ref{F:Welded}. Performing a $(F1)$ move on~$\beta$, one obtains a strand diagram $\beta^\prime$ whose associated braid Gauss diagram is~$g^\prime$. 
\endproof
%

\section{An historical note and other references}
\label{S:other}

Loop braids have first been considered as \emph{motions} of $\nn$ unknotted, unlinked circles in~$\Rr^3$, as seen in Section~\ref{S:MCG}, a notion that can be translated in terms of mapping classes. More in detail, Fox and Neuwirth give in 1962 a new proof of the standard presentation of the Artin braid groups $\BB\nn$, in terms of fundamental groups of the spaces of configurations of $\nn$ unordered points in the plane. From this, the interpretation of braids as \emph{motions} on $\nn$ points on a plane arises. The same year Dahm, in his unpublished Ph.D thesis~\cite{Dahm}, generalises this concept defining the group of motions of a compact subspace $N$ in a manifold $M$. He  applies this construction to a collection of $\nn$ unknotted, unlinked circles in the $3$-dimensional space,  and proves the resulting mapping class groups to be isomorphic to particular subgroups of the groups of automorphisms of the free group~$\Aut(\F\nn)$. This is the first appearance of the loop braid groups, and his results are published and extended by Goldsmith~\cite{GOLD}. Some insight on the link between motion groups and mapping class groups can be found in Bellingeri-Cattabriga's work~\cite[Section~5]{BC}.

In 1986 McCool~\cite{McC} considers the subgroups of \emph{basis-conjugating automorphisms} of the groups $\Aut(\F\nn)$ and gives a presentation for them. We have seen in Section~\ref{S:automorphisms} that these are isomorphic to the pure loop braid groups~$\PLB\nn$. A decade later Savushkina~\cite{SAV} studies the groups of \emph{permutation-conjugacy automorphisms}, isomorphic to the loop braid groups~$\LB\nn$. She gives a presentation for these and investigates many properties of these groups.

Loops braids are introduced then as \emph{welded braids}, as defined by Fenn-Rim\'{a}nyi-Rourke~\cite{FRR}, in the form of a presented group, with diagrams representing its generators, and ten years later Baez-Crans-Wise~\cite{BWC} establish an isomorphism between the motion groups of $\nn$ non-flipping circles and the groups of welded braids.

More recently, in~\cite{BH}, we find an interpretation of these groups in terms of fundamental groups of the configuration spaces of $\nn$ unlinked Euclidean circles. About the topological interpretation of loop braids in terms of braided ribbon tubes, one can find details in~\cite{Kam2}: this interpretation is also widely used in~\cite{ABMW}. Finally, as explained in~\cite{Bn}, they have also an interpretation in terms of Gauss braid diagrams: a point of view that allows us to consider the loop braid groups as quotients of the virtual braid groups.

Several computations and conjectures about the cohomology algebras of the pure loop braid groups appear in the literature. Investigating the properties of these algebra, for example in relation with resonance varieties and lower central series ranks, allows us to compare the loop braid groups to the braid groups, and to other generalizations of the braid groups. About this topic we refer to~\cite{JCM, CPVW, SW, Ver02}.

A word should be said about the theory of representations of the loop braid groups, which appears to be an upcoming topic. Burau representation extends trivially to loop braid groups using Magnus expansion and Fox derivatives~\cite{Bar}, however it is still unknown if the loop braid groups are linear. Some new results on local representations of loop braid groups, rising as extensions of braid groups representations, can be found in~\cite{KMRW} and~\cite{BCH}. However the study of finite dimensional quotients of algebras of loop braid groups is yet to be found in the literature. In~\cite{KMRW} the authors show interest also in certain remarkable quotients of loop braid groups, the \emph{symmetric loop braid groups} (also known as \emph{unrestricted virtual braid groups}~\cite{KL}).  These groups have a very simple algebraic structure and we refer to \cite{BBD,Nas, ABMW3} for their properties and applications to \emph{fused links}.

\section*{Acknowledgements}
I would like to thank P.~Bellingeri for his invaluable advice and insight throughout this work. I would also like to thank: B.~Audoux for accurate comments on the manuscript; A.~Hatcher for useful remarks and for pointing out precious references; S.~Kamada and D.~Silver for the comments on my thesis, which contains an expanded version of this paper; and T. Brendle, J.-B.~Meilhan, and E.~Wagner for helpful conversations and comments.

\bibliography{JourneyLoopBraid.bib}{}
\bibliographystyle{plain}

\end{document}